\documentclass[12pt,hidelinks]{article}
\usepackage[utf8]{inputenc}
\usepackage[T1]{fontenc}
\usepackage{amsmath, amssymb, amsfonts, url, float, graphicx, float}
\usepackage{fullpage}
\usepackage{xcolor}
\usepackage[small,bf]{caption}
\usepackage{cite}
\usepackage{tikz, pgfplots}
\pgfplotsset{
    compat=1.16,
}
\usetikzlibrary{topaths,intersections,calc}
\usetikzlibrary{positioning,chains,fit,shapes,calc}
\usepackage{graphicx}
\graphicspath{{figures/}}
\usepgfplotslibrary{fillbetween}
\newcommand{\arb}{\textbf{arb}}

\newcommand{\downto}{\downarrow}

\let\emptyset\varnothing
\let\phi\varphi

\newcommand{\ones}{\mathbf 1}
\newcommand{\reals}{\mathbf{R}}

\newcommand{\cl}{{\mathop{\bf cl}}} 

\newcommand{\intr}{\mathop{\bf int}}

\newcommand{\cf}{{\it cf.}}
\newcommand{\eg}{{\it e.g.}}
\newcommand{\ie}{{\it i.e.}}

\newcommand{\BEAS}{\begin{eqnarray*}}
\newcommand{\EEAS}{\end{eqnarray*}}
\newcommand{\BEA}{\begin{eqnarray}}
\newcommand{\EEA}{\end{eqnarray}}
\newcommand{\BEQ}{\begin{equation}}
\newcommand{\EEQ}{\end{equation}}
\newcommand{\BIT}{\begin{itemize}}
\newcommand{\EIT}{\end{itemize}}

\title{The Geometry of Constant Function Market Makers}
\author{
    Guillermo Angeris\thanks{The authors are listed in alphabetical order.}\\\texttt{\small gangeris@baincapital.com} \and
    Tarun Chitra\\\texttt{\small tarun@gauntlet.network} \and
    Theo Diamandis\\\texttt{\small tdiamand@mit.edu} \and
    Alex Evans\\\texttt{\small aevans@baincapital.com} \and
    Kshitij Kulkarni\\\texttt{\small ksk@eecs.berkeley.edu}
}
\date{July 2023}

\begin{document} 
\maketitle 

\begin{abstract}
    Constant function market makers (CFMMs) are the most popular type of
    decentralized trading venue for cryptocurrency tokens. In this paper, we
    give a very general geometric framework (or `axioms') which encompass and
    generalize many of the known results for CFMMs in the literature, without
    requiring strong conditions such as differentiability or homogeneity. One
    particular consequence of this framework is that every CFMM has a (unique)
    canonical trading function that is nondecreasing, concave, and homogeneous,
    showing that many results known only for homogeneous trading functions are
    actually fully general. We also show that CFMMs satisfy a number of
    intuitive and geometric composition rules, and give a new proof, via conic
    duality, of the equivalence of the portfolio value function and the trading
    function. Many results are extended to the general setting where the CFMM
    is not assumed to be path-independent, but only one trade is allowed.
    Finally, we show that all `path-independent' CFMMs have a simple geometric
    description that does not depend on any notion of a `trading history'.
\end{abstract}

\section*{Introduction}
The study of automated market makers has existed for many decades, with
roots in the scoring rule literature dating back to at least the 
1950s~\cite{mccarthy1956measures}. However, these mechanisms only reached
mass adoption after being implemented as decentralized exchanges on
blockchains. Surprisingly, the types of automated market maker that are most
popular in practice bear little resemblance to those proposed prior to the
invention of blockchains. Instead, the most popular blockchain-based automated
market makers are what are known as the \emph{constant function market makers},
or CFMMs. These market makers are generally simpler than earlier market maker
designs, such as those based on the logarithmic scoring rule, and provide a
means for efficient liquidity aggregation and order routing. But why have these
mechanisms succeeded?

One of the main reasons that these mechanisms have been so popular in the
cryptocurrency space is that solving the optimal arbitrage problem---the
problem of how much to trade in order to equalize prices between CFMMs and
other venues---is generally (computationally) `easy'~\cite{angeris2022optimal}.
This ease comes directly from the fact that CFMMs satisfy a general, very
geometric, notion of convexity. Though the initial line of work, which defined
CFMMs as a useful class, focused on their geometric
properties~\cite{angeris2020improved}, the majority of research on CFMMs has
focused on analytic properties of CFMMs that depend on explicit
parameterizations~\cite{lehar2021decentralized, wu2022constant,
milionis2023complexityapproximation, schlegel2023axioms,
frongillo2023axiomatic, milionis2023myersonian, goyal2023finding}.

There are important reasons to examine geometry of CFMMs directly. First, a
geometric lens leads to very natural statements for many of the properties of,
and operations one can perform on, CFMMs. Second, many `surprising' decisions
made by developers that `worked in practice' can be explained by understanding
the geometry of CFMMs; for example, intuitively, the `curvature' of a CFMM
corresponds to a notion of liquidity~\cite{angeris2022when}, which was known by
practitioners well before its formalization. Third, the geometric setting for
CFMMs is very general and rarely requires notions of differentiability,
homogeneity, or other similar properties. Finally, the geometric view of CFMMs
allows for reasoning about CFMMs without regards to its particular trading
function and/or its representation.

\paragraph{What is a CFMM?} CFMMs are relatively simple to describe
analytically, which is likely why many researchers and practitioners work with
them in this `analytic' setting. A CFMM consists of two main objects: a trading
function $\phi: \reals_+^n \to \reals$ (sometimes called an `invariant') and a
vector of reserves $R \in \reals_+^n$. A user proposes a trade, represented as
a portfolio $\Delta \in \reals^n$, and the trade is valid if the trading
function, evaluated on the reserves after the trade is completed, has the same
value as the function evaluated on the reserves before the trade is completed,
\ie, if $\phi(R - \Delta) = \phi(R)$. (Hence the name `constant function market
maker'.) If this equality holds, the CFMM then pays out $\Delta$ to the user and
resulting in new reserves $R - \Delta$. (If the trade is invalid, nothing is
paid out or received from the user.) Liquidity providers, who provide the
reserves $R$ against which trades are made, earn fees from these trades. The
fact that this process is simple to describe and implement, along with having
many strong theoretical guarantees, has been part of its reason for success,
especially within difficult-to-secure environments such as public blockchains.
Despite their simple description, CFMMs have spawned a large amount of research
into their financial, arbitrage, and routing properties
(\eg,~\cite{angeris2020improved, danos2021global, diamandis2023efficient,
fukasawa2023weighted, mazorra2023towards}, among many others).

\paragraph{Analytic vs.\ geometric properties.} Many descriptions of CFMMs
use a coordinate-dependent (or `analytic') version of CFMMs, focusing on
the representation of $\phi$ as an explicit function. For instance, Uniswap is
commonly described via the trading function $\phi(R) = \prod_{i=1}^n R_i$.
However, there are a number of equivalent trading functions such as $\tilde \phi(R) =
\left(\prod_{i=1}^n R_i\right)^c$ for any $c > 0$. While the pricing and
behavior is equivalent for these different representations, the mechanics of
many definitions hinge upon the specific representation provided. (Some might
demand concavity or monotonicity, for example.) We call any properties that are
dependent on the particular representation of the trading function
\emph{analytic} properties.

On the other hand, nearly by definition, geometric descriptions of CFMMs are
unique and relatively simple to handle. For instance, there is a natural
`addition' operator for CFMMs using a geometric representation. Describing the
corresponding operation on trading functions is not obvious and likely has no
natural analogue. (We do show that there is another functional representation,
given by the portfolio value function, that does have a natural
correspondence.) This idea that certain operations such as addition, are `easy'
to perform on CFMMs, when defined geometrically, is one of the reasons that
proofs using geometry can be significantly more succinct than those using
analytic means.

\paragraph{This paper.} In this paper, we focus on representing CFMMs via
classical geometric objects such as convex sets and cones, assuming a bare
minimum of requirements. Using these objects, we replicate the results of a
number of papers for CFMMs without fees (also known as the `path-independent'
CFMMs) and many results in the case of a single trade with no restrictions on
the CFMM. The key objects we look at are particular cones which we call the
liquidity cones, and their corresponding conic duals. We construct many of the
`usual objects' such as trading functions, portfolio value functions, no
arbitrage intervals, and so on, directly from these objects. This leads to a
number of interesting results: for example, that every (path-independent) CFMM
has a canonical trading function that is nondecreasing, concave, and
homogeneous, along with new proofs for older, previously known results. We
assume a reasonable amount of familiarity with convex optimization and provide
a very short primer on conic duality in appendix~\ref{app:conic-duality} as a
refresher.

\section{Fee-free constant function market makers}
In this section, we consider the general case of constant function market
makers that are path-independent. We show the connection between these
`path-independent' or `fee-free' constant function market makers and `general'
constant function market makers later, in~\S\ref{sec:path-independence}. We
consider this case first as this is the most common case in the
literature~\cite{angeris2020improved,frongillo2023axiomatic,goyal2023finding},
and is a good starting point for the more general case.

\paragraph{Section layout.} The section begins with a basic set of requirements
(sometimes called `axioms') which are of a different form than the standard
assumptions made in many texts. We will show that from these requirements,
which are mostly geometric in origin, we can derive many known results and a
number of generalizations that, to our knowledge, are not known in the overall
literature. For example, one important case is that any CFMM has a canonical
trading function that is homogeneous, nondecreasing, and concave. This is
usually taken as an assumption in some form (see,
\eg,~\cite{angerisNotePrivacyConstant2021, frongillo2023axiomatic,
schlegel2023axioms}), but we show here that it is true of any CFMM satisfying
some basic properties that are essentially necessary for a CFMM to be
reasonable. (Indeed, these properties are almost always part of a much longer
list, or are easy consequences of a subset of assumptions generally made in the
literature.) This geometric set up also simplifies a number of known statements
in the literature, such as those of~\cite{angeris2023replicating}, by showing
that the equivalence of a portfolio value function and a trading function is a
special case of conic duality.

\subsection{Reachable set}\label{sec:reachable-set}
We will define the reachable set of reserves as a set $S\subseteq \reals^n$
satisfying certain requirements. This set will represent the valid holdings of
a constant function market maker (CFMM). In general, if $R \in S$ are the
current reserves of the constant function market maker, then any trader may
change the reserves to $R' \in S$ by selling $R' - R$ to the CFMM. The trader
would then receive the entries of $R' - R$ which are negative, and tender the
entries which are positive to the CFMM. In a certain sense, we may view the
reachable set $S$ as the set of valid states available to the CFMM.

\paragraph{Definition.} We say a set $S$ is a \emph{reachable set} (which
defines a fee-free, or `path independent' CFMM) if it satisfies these rules or `axioms':
\begin{enumerate}
    \item All reserves are nonnegative; that is, $S \subseteq \reals_+^n$
    \item The set $S$ is nonempty, closed, and convex
    \item The set $S$ is \emph{upward closed}; \ie, if $R \in S$, then any $R' \ge R$ has $R' \in S$
\end{enumerate}
From these three rules, we will recover (and generalize) many of the results
known in the literature. In general, while we do not assume that $0 \not \in S$,
we note that this is a silly case as we would then have $S = \reals_+^n$,
so this case is often excluded from many of the proofs presented.

\paragraph{High-level interpretation.} The first requirement means that a
constant function market maker cannot take on debt, or that the position is
always solvent. Many, but not all, results hold with some slight modifications,
even in the case where this condition is relaxed. The convexity requirement
roughly corresponds to the fact that increasing the size of a trade does not
result in a better exchange rate for the trader. The nonemptyness of $S$ just
means that $S$ is nontrivial, while the closedness is a technical condition.
Finally, the `upwards closed' condition means that, if a CFMM accepts some
trade, then it would always accept a different trade that tenders more of any
asset. (This condition is not technically necessary: it suffices that, given a
nonempty set $S$ satisfying the first condition, the set $\tilde S = S +
\reals_+^n$ satisfies the second condition above. Almost all results shown
below hold in this case.) The last condition also lets us interpret the
boundary of $S$ as a Pareto-optimal frontier for the possible reserves in the
sense that no rational trader would ever trade on the interior of $S$.

\begin{figure}
    \centering
    \begin{tikzpicture}[scale=1.0]
    \begin{axis}[
        axis on top=true,
        xmin=-1, xmax=4,
        ymin=-1, ymax=4,
        axis lines=center,
        x label style={anchor=north},
        y label style={anchor=east},
        xlabel={$R_1$},
        ylabel={$R_2$},
        grid=major,
        legend pos=outer north east,
        samples=200,
        domain=-4:4,
        xtick=\empty,
        ytick=\empty
        ]

        \addplot[black, very thick, domain=0.01:4.99] {1/x};
        \addplot[gray!30, domain=0.01:4.99, fill, opacity=0.5, draw=none] (1/x, x) -| (4,4) -- cycle;
        \node at (2.25, 2.25) {\tiny $S = \{R \in \reals_+^2 \mid R_1R_2 \ge 1\} $};
    \end{axis}
\end{tikzpicture}
    \hfill
    \begin{tikzpicture}[scale=1.0]
    \begin{axis}[
        axis on top=true,
        xmin=-1, xmax=4,
        ymin=-1, ymax=4,
        axis lines=center,
        x label style={anchor=north},
        y label style={anchor=east},
        xlabel={$R_1$},
        ylabel={$R_2$},
        grid=major,
        legend pos=outer north east,
        samples=200,
        domain=-4:4,
        xtick=\empty,
        ytick=\empty
        ]

        \addplot[black, very thick, domain=0.0:1.5] {1/(x + 1/2) - 1/2};
        \addplot[black, very thick, domain=1.5:4.0] {0};
        \draw[black, very thick] {(0, 1.5) -- (0, 4)};
        \addplot[gray!30, domain=0.0:1.5, fill, opacity=0.5, draw=none] (1/(x + 1/2) - 1/2, x) -| (0,4) -| (4,4) -| (4,0) -- cycle;
        \node at (2.03, 2.0) {\tiny $S = \{R \in \reals_+^2 \mid (R_1 + \tfrac{1}{2})(R_2 + \tfrac{1}{2}) \ge 1\} $};
    \end{axis}
\end{tikzpicture}
    \caption{The set of reachable reserves for Uniswap (left) and Uniswap v3 (right).}
    \label{fig:uni-reserves}
\end{figure}
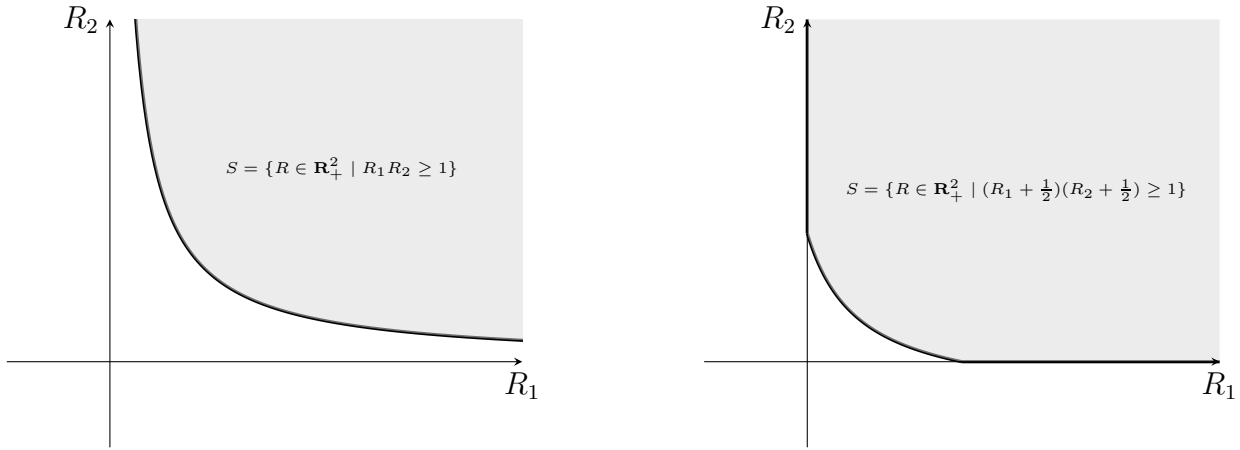

\paragraph{Examples.} One of the canonical examples of a reachable set is that
of Uniswap~\cite{adams2020uniswap}, defined
\[
    S = \{R \in \reals_+^2 \mid R_1R_2 \ge k\},
\]
where $k > 0$ is a constant. See figure~\ref{fig:uni-reserves} for an example.
Another example is that of a `tick' in Uniswap v3~\cite{adams2021uniswap},
which is defined
\[
    S = \{R \in \reals_+^2 \mid (R_1 + \alpha)(R_2 + \beta) \ge k\},
\]
where, again, $\alpha, \beta, k > 0$ are some provided constants.

\paragraph{Quasiconcavity.} Note that, in these examples, $S$ is the superlevel
set of some quasiconcave, nondecreasing function. In fact, we can show that any
nonempty set $S$ defined by
\begin{equation}\label{eq:quasiconcave}
    S = \{R \in \reals_+^n\mid \psi(R) \ge \alpha\},
\end{equation}
with quasiconcave, nondecreasing $\psi: \reals_+^n \to \reals \cup
\{-\infty\}$, generates a reachable set satisfying the required conditions.
This includes~\cite{schlegel2023axioms} and~\cite{frongillo2023axiomatic} as a
special case, though we do not require homogeneity. (Indeed, homogenity is not
needed as an assumption as we will later show that one can always choose $\psi$
to be concave, nondecreasing, and homogeneous for any set $S$ satisfying the
reachable set conditions, even when the `original' function $\psi$ is not.) We
may also replace the inequality with an equality and define the set
\[
    S = \{R' \in \reals_+^n \mid R' \ge R, ~\text{for some} ~\psi(R) = \alpha\}.
\]
Note that these two definitions are equivalent if $\psi$ is continuous in some
neighborhood $\psi^{-1}(N)$ where $N$ is a neighborhood around $\alpha$.

\subsection{Composition rules}\label{sec:composition-rules}
An interesting consequence of the definition of reachable sets is that
these sets, and therefore CFMMs, satisfy certain composition rules, some of
which were known in the literature under additional assumptions~\cite{engel2021}.
These rules follow directly from the calculus of convex
sets~\cite[\S2.3]{cvxbook} and require no additional assumptions than those
given in~\S\ref{sec:reachable-set}.

\paragraph{Nonnegative scaling.} Given a reachable set $S$, we may scale the
set by $\alpha \ge 0$ to get $\alpha S$, which is another reasonable reachable
set satisfying the conditions. (We will see later
in~\S\ref{sec:liquidity-provision} that this scaling corresponds to adding or
removing liquidity to the CFMM.)

\begin{figure}
    \centering
    \begin{tikzpicture}[scale=0.7]
    \begin{axis}[
        axis on top=true,
        xmin=-1, xmax=3,
        ymin=-1, ymax=3,
        axis lines=center,
        x label style={anchor=north},
        y label style={anchor=east},
        xlabel={$R_1$},
        ylabel={$R_2$},
        grid=major,
        legend pos=outer north east,
        samples=200,
        domain=-4:4,
        xtick=\empty,
        ytick=\empty
        ]

        \addplot[black, very thick, domain=0.0:0.5] {1/(x + 1/2) - 1};
        \addplot[black, very thick, domain=0.5:4.0] {0};
        \draw[black, very thick] {(0, 1.0) -- (0, 4)};
        \addplot[gray!30, domain=0.0:1.0, fill, opacity=0.5, draw=none] (1/(x + 1) - 1/2, x) -| (0,4) -| (4,4) -| (4,0) -- cycle;
        \node at (1.5, 1.5) {$S_1$};
    \end{axis}
\end{tikzpicture}
\hfill
\begin{tikzpicture}[scale=0.7]
    \begin{axis}[
        axis on top=true,
        xmin=-1, xmax=3,
        ymin=-1, ymax=3,
        axis lines=center,
        x label style={anchor=north},
        y label style={anchor=east},
        xlabel={$R_1$},
        ylabel={$R_2$},
        grid=major,
        legend pos=outer north east,
        samples=200,
        domain=-4:4,
        xtick=\empty,
        ytick=\empty
        ]

        \addplot[black, very thick, domain=0.0:1.0] {1/(x + 1) - 1/2};
        \addplot[black, very thick, domain=1.0:4.0] {0};
        \draw[black, very thick] {(0, 0.5) -- (0, 4)};
        \addplot[gray!30, domain=0.0:0.5, fill, opacity=0.5, draw=none] (1/(x + 1/2) - 1, x) -| (0,4) -| (4,4) -| (4,0) -- cycle;
        \node at (1.5, 1.5) {$S_2$};
    \end{axis}
\end{tikzpicture}
\hfill
\begin{tikzpicture}[scale=0.7]
    \begin{axis}[
        axis on top=true,
        xmin=-1, xmax=3,
        ymin=-1, ymax=3,
        axis lines=center,
        x label style={anchor=north},
        y label style={anchor=east},
        xlabel={$R_1$},
        ylabel={$R_2$},
        grid=major,
        legend pos=outer north east,
        samples=200,
        domain=-4:4,
        xtick=\empty,
        ytick=\empty
        ]

        \addplot[black, very thick, domain=0.0:1.5] {1/(x + 1/2) - 1/2};
        \addplot[black, very thick, domain=1.5:4.0] {0};
        \draw[black, very thick] {(0, 1.5) -- (0, 4)};
        \addplot[gray!30, domain=0.0:1.5, fill, opacity=0.5, draw=none] (1/(x + 1/2) - 1/2, x) -| (0,4) -| (4,4) -| (4,0) -- cycle;
        \node at (1.5, 1.5) {$S_1 + S_2$};
    \end{axis}
\end{tikzpicture}
    \caption{Adding two Uniswap v3 bounded liquidity pools (left, middle) gives
    us another CFMM (right).}
    \label{fig:minkowski-reserves-addition}
\end{figure}
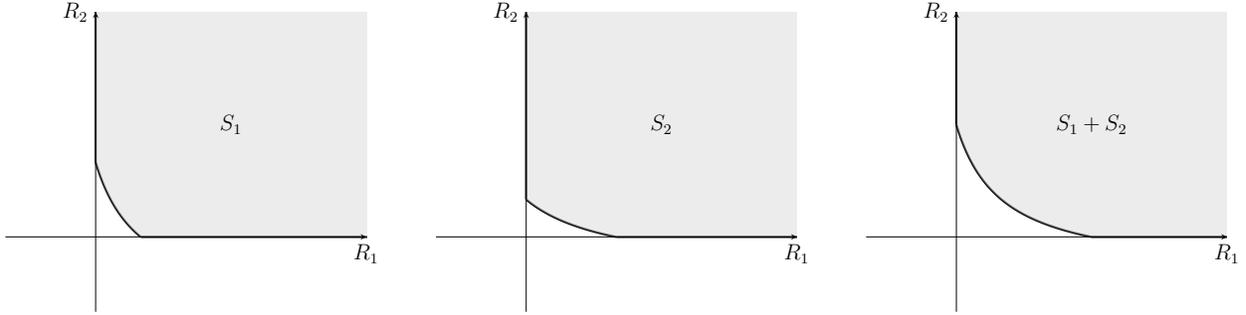

\paragraph{Set addition.} We may also add any two reachable sets $S$ and $S'$,
which gives another reachable set
\[
    S + S' = \{R + R' \mid R \in S, ~ R' \in S'\}.
\]
This set is convex and nonempty, and it is not hard to prove the set is closed
since $S$ and $S'$ are both contained in the positive orthant. It is also clear
that $S + S'$ is upward closed since each of $S$ and $S'$ are upward closed.
These sums have the `simple' interpretation that $S+S'$ are the possible
combined holdings of the two CFMMs. Additionally, this, combined with
nonnegative scaling, means that taking nonnegative linear combinations of
trading sets always yields another trading set. We provide an example in
figure~\ref{fig:minkowski-reserves-addition}.

\paragraph{Nonnegative matrix multiplication.} Another important rule is that
multiplication by a nonnegative matrix $A \in \reals_+^{n \times p}$ and
`upwards closure' of the resulting set gives another reachable set; \ie, the
set
\[
    AS + \reals_+^n = \{R' \in \reals_+^n \mid R' \ge AR ~ \text{for some} ~ R \in S\},
\]
is a reachable set. This operation can be interpreted when looking at each row
$j=1, \dots, n$ of $A$, which we write as $\tilde a_j^T$. Given some vector $R
\in S$, then $\tilde a_j^TR = (AR)_j$. This entry, $(AR)_j$, can then be seen
as a type of `meta-asset', whose value is equal to a weighted basket of assets,
where the weights are the entries of $\tilde a_j$. This is a reachable set
since $AR \in \reals^n_+$ for any $R \in S \subseteq \reals^n_+$ and $AS$
is a convex set if $S$ is convex. (The set is clearly upward closed by definition.)

\paragraph{Special case: projection.} An important special case is when the
matrix $A$ projects all components of a trading set into a larger space.
More specifically, let $A$ be a matrix of the form
\[
    A = \begin{bmatrix}
        a_1 & a_2 & \cdots & a_k
    \end{bmatrix}
\]
with $a_i \in \reals^n$ and $k \le n$ being all distinct unit basis vectors
(\ie, $a_i$ is $0$ everywhere except at exactly one entry, where it is $1$). We
can interpret $AS + \reals_+^k$ in the following way: if there is a `list' of
$n$ assets, and the CFMM defined by $S$ trades only $k$ of those assets, then
$AS$ is a CFMM which trades these $k$ assets and zero of the remaining possible
$n-k$ assets. (The CFMM will happily accept any of the remaining $n-k$ assets,
but tender nothing for them: a trade no rational user would want.)

\paragraph{Intersection.} Finally, we can take the intersection of reachable
sets, which yields another reachable set; \ie, if $S$ and $S'$ are reachable
sets then $S\cap S'$ is similarly a reachable set. This corresponds to a CFMM
whose reachable reserves can only be those which the individual CFMMs have in
common. Though this is not a natural operation for CFMMs which already exist on
chain, it is a useful theoretical operation for constructing CFMMs with
particular properties. (Indeed, we will see that constructing a CFMM from a
portfolio value function, to be presented later, is possible only due to this
intersection property.)

\paragraph{Aggregate CFMMs.} Combining the previous two rules gives us a very
general way of `combining' CFMMs which trade different (but potentially
overlapping) baskets of assets. Assume we have $m$ constant function market
makers and a universe of $n$ assets. We will have CFMMs $i=1, \dots, m$ with
reachable sets $S_i \subseteq \reals_+^{n_i}$, each trading a subset tokens of
$n_i$ tokens. We introduce matrices $A_i\in \reals_+^{n \times n_i}$ which 
map the `local' basket of $n_i$ tokens for CFMM $i$ to the global universe of 
$n$ tokens. We have $(A_i)_{jk} = 1$ if token $k$ in market $i$'s local index
corresponds to global token index $j$, and $(A_i)_{jk}=0$ otherwise. We note
that the ordering of tokens in the local index does not need to be the same as
the global ordering. Then the set
\[
    \tilde S = \sum_{i=1}^m A_i S_i
\]
is a \emph{aggregate CFMM} which corresponds exactly to the set of all possible
holdings for every CFMM in the network. Such CFMMs were first implicitly
defined for Uniswap v3~\cite{adams2021uniswap}, and later used
in~\cite{chitra2021liveness} to prove some basic approximation bounds,
while~\cite{milionis2023complexityapproximation} defined a notion of
`complexity' based on similar ideas, and, finally~\cite{diamandis2023efficient}
defined them as part of the solution method for optimal routing.

\paragraph{Extensions to negative reserves.} There are some basic
generalizations of some of these conditions in the case where the set $S$ is
not contained in the positive orthant. In this case, the CFMM can take on debt.
If the debt is unbounded, it is possible to create sets $S$ and $S'$ such that
$S+S'$ is not closed, so the resulting set would not be a reachable set. On the
other hand, it is not hard to show that allowing bounded debt (\ie, there
exists some $x \in \reals_+^n$ such that $x + S \subseteq \reals_+^n$) means
that an analogous statement does still hold by a nearly identical proof.

\subsection{Liquidity cone and canonical trading function}
In this subsection we introduce the liquidity cone for a reachable set $S$. The
liquidity cone is a kind of `homogenized' version of the reachable set defined
previously that simplifies a number of later derivations. Its definition will
also suggest a canonical trading function: a trading function that corresponds
to the reachable set $S$ and is nondecreasing, homogeneous, and concave.

\subsubsection{Liquidity cone} 
The \emph{liquidity cone} for reachable set $S$ is defined as
\begin{equation}\label{eq:liquidity-cone}
    K = \cl\{(R, \lambda) \in \reals^{n+1} \mid R/\lambda \in S, ~ \lambda > 0\},
\end{equation}
where $\cl$ is the closure of the set. The set $K$ is a cone as $(R, \lambda)
\in K$ implies that $(\alpha R, \alpha \lambda) \in K$ for any $\alpha \ge 0$.
The name `liquidity cone' comes from the fact that, if $(R, \lambda) \in K$
then the largest such $\lambda$ indicates, roughly speaking, the amount of
liquidity available from reserves $R$. (We will see what this means in a later
section.)

\paragraph{Basic properties.} The liquidity cone $K$ has some important
properties we use later in this section. First, the set $K$ is nonempty as $S$
is nonempty and $S\times \{1\} \subseteq K$. We also have that $0 \in K$ as $K$
is nonempty and closed. To see this, if $y\in K$ then $\alpha y \in K$, so
$\alpha \downto 0$ gives the result. The cone $K$ is also a convex cone as it
is the closure of the perspective transform on the convex set $S\times
\reals_{++}$ (see, \eg,~\cite[\S 2.3.3]{cvxbook}).

\paragraph{Upward closedness.} The cone $K$ is not upward closed, but is
`almost upward closed' in the following sense: if $(R, \lambda) \in K$ and $R'
\ge R$ with $\lambda' \le \lambda$ then $(R', \lambda') \in K$. In particular,
note that the inequality over $\lambda$ is reversed. Showing this fact is just
a definitional exercise.

\paragraph{Positive reachability.} We also have that,
\begin{equation}\label{eq:positive-reachability}
    (\reals_{++}^n, 0) \subseteq K.
\end{equation}
This follows from the fact that the set $S$ is nonempty. To see this, let $R
\in S$ and note that, for any strictly positive vector $R' \in \reals_{++}^n$
we know that $R'/\lambda \ge R$ for $\lambda$ small enough, so $(R', \lambda)
\in K$. Finally, since $(R', \lambda) \in K$ implies that $(R', \lambda') \in
K$ for any $\lambda' \le \lambda$, then we are done by setting $\lambda' = 0$.
Roughly speaking, this corresponds to the intuitive fact that every nonnegative
basket is a feasible set of reserves, at some `large enough' multiple. 
This observation is taken as an assumption in~\cite{frongillo2023axiomatic}
and~\cite{schlegel2023axioms}, but is a direct consequence of the definition of
the reachable set. Additionally, since $K$ is closed we have
\[
    (\reals_{+}^n, 0) \subseteq K,
\]
though this construction is less useful than the previous.

\paragraph{Reachable set.} We may, of course recover the reachable set
from the liquidity cone in a variety of ways. Perhaps the simplest is to
note that, for any $\lambda > 0$ we have
\begin{equation}\label{eq:reachable-set-cone}
    S = \{R/\lambda \mid (R, \lambda) \in K\}.
\end{equation}
This is easy to see as $(R, \lambda) = \lambda(R/\lambda, 1) \in K$, and,
since $K$ is a cone, this is if, and only if, $(R/\lambda, 1) \in K$ which is
also if, and only if, $R/\lambda \in S$. This will be useful in what follows.

\subsubsection{Canonical trading function}\label{sec:trading-function}
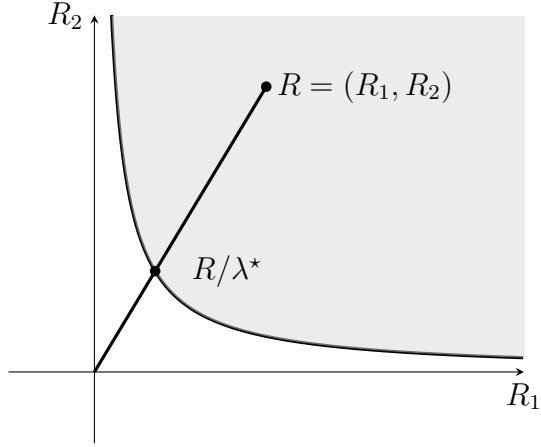
\begin{figure}
    \centering
    \begin{tikzpicture}[scale=1.0]
    \begin{axis}[
        axis on top=true,
        xmin=-1, xmax=5,
        ymin=-1, ymax=5,
        axis lines=center,
        x label style={anchor=north},
        y label style={anchor=east},
        xlabel={$R_1$},
        ylabel={$R_2$},
        grid=major,
        legend pos=outer north east,
        samples=200,
        domain=-5:5,
        xtick=\empty,
        ytick=\empty
        ]

        \addplot[black, very thick, domain=0.01:4.99] {1/x};
        \addplot[gray!30, domain=0.01:4.99, fill, opacity=0.5, draw=none] (1/x, x) -| (5,5) -- cycle;
        \addplot[black, very thick, domain=0.0:2.0] {2*x};
        \node[above, right] at (2, 4) {$R = (R_1, R_2)$};
        \fill (2,4) circle[radius=2pt];


        \node[right] at (0.7071, 1.4142) {\hspace{0.5em} $R / \lambda^\star$};
        \fill (0.7071, 1.4142) circle[radius=2pt];
    \end{axis}
\end{tikzpicture}
    \caption{Another interpretation of the canonical trading
        function~\eqref{eq:canonical-cone}: we scale along the line segment defined by
        $(R_1, R_2)$ to $(0,0)$, with scale factor $1/\lambda$, increasing $\lambda$
        until we hit the reachable set boundary.}
    \label{fig:uni-trading-function}
\end{figure}

Given any liquidity cone $K$ for a reachable set $S$, we will define a
\emph{canonical trading function},
\begin{equation}\label{eq:canonical-cone}
    \phi(R) = \sup\{\lambda \mid (R, \lambda) \in K\},
\end{equation}
setting $\phi(R) = 0$ if the set is empty. (Since $K$ is closed, we may replace
the $\sup$ with a $\max$ if $0 \not \in S$, which we assume for the remainder
of the section.) In terms of the trading set $S$, we may write this as
\[
    \phi(R) = \sup\{\lambda > 0 \mid R/\lambda \in S\},
\]
using the definition of the liquidity cone $K$. If the reachable
set $S$ is written using a nondecreasing, quasiconcave, but not necessarily
concave, function as in~\eqref{eq:quasiconcave}, then we can `canonicalize'
this trading function by writing
\begin{equation}\label{eq:canonical-reach}
    \phi(R) = \sup\{\lambda > 0 \mid \psi(R/\lambda) \ge k\}.
\end{equation}
Note that, if $\psi$ is continuous, this is the same as finding the largest
positive root over $\lambda$ of $\psi(R/\lambda) = k$. If the function is
strictly increasing (as is often the case) then the positive root is unique and
it suffices only to find it. Figure~\ref{fig:uni-trading-function} illustrates
this definition for the case of Uniswap.

\paragraph{Computational considerations.} It may be the case that the canonical
trading function~\eqref{eq:canonical-reach} has no closed form solution. From
the previous, since we know that computing the value of the canonical trading
function at some reserves $R$ corresponds to a root-finding problem, we may do
this using efficiently by using bisection (as $\psi$ is assumed to be
nondecreasing) or, if $\psi$ is differentiable, using Newton's method for
finding the positive root. In either case, computing $\phi(R)$ can be done
efficiently in practice. (As a side note: if bisection is used, it suffices
to run it only until the bracketing interval is either fully contained in $[0, 1)$
or $[1, \infty)$. In the former, the reserves are guaranteed to be infeasible, while
in the latter they are guaranteed to be feasible.)

\paragraph{Reachable set.}
From~\eqref{eq:reachable-set-cone} we can recover the set $S$ from
this canonical trading function since
\[
    S = \{R \in \reals_+^n \mid \phi(R) \ge 1\}.
\]
(Of course, the set of $R$ such that $\phi(R) = 1$ gives the boundary of $S$.)
Additionally, note that if $\phi(R) > 0$, which is always true if $R \in
\reals_{++}^n$ is strictly positive, from positive
reachability~\eqref{eq:positive-reachability}, then
\begin{equation}\label{eq:reachable-scale}
    \frac{R}{\phi(R)} \in S.
\end{equation}

\paragraph{Concavity.} This function is concave, as it is the partial
maximization of the concave function
\[
    f(R, \lambda) = \lambda - I(R, \lambda)
\]
over $\lambda$, where $I$ is the indicator function of the (convex) set $K$, defined
$I(R, \lambda) = 0$ if $(R, \lambda) \in K$ and $+\infty$ otherwise.

\paragraph{Homogeneity.} The trading function $\phi$ is
homogeneous for $\alpha > 0$ since
\[
    \phi(\alpha R) = \sup\{\lambda \mid (\alpha R, \lambda) \in K\}.
\]
Since $K$ is a cone, then $(\alpha R, \lambda) \in K$ if, and only if, $(R, \lambda/\alpha) \in K$.
Setting $\bar \lambda = \lambda/\alpha$, then we have
\[
    \phi(\alpha R) = \sup\{\alpha \bar \lambda \mid (R, \bar \lambda) \in K\} = \alpha \phi(R).
\]
For $\alpha = 0$ the result follows since $0 \in K$.

\paragraph{Monotonicity.} The trading function is nondecreasing from the
`almost upward closed' property mentioned previously. For the remainder of the
paper, we will call a function that is concave, homogeneous, and nondecreasing
a \emph{consistent} function.

\paragraph{Marginal prices.}\label{sec:marg-prices}
Given $R$ with $\phi(R) = 1$, \ie, the starting reserves are `reasonable' and
$\phi$ differentiable at $R$, then, from concavity,
\[
    \phi(R + \Delta) \le \phi(R) + \nabla\phi(R)^T\Delta.
\]
(We may replace the gradient with a supergradient for a more general condition.)
If the trade $\Delta$ is feasible in that $\phi(R + \Delta) \ge \phi(R)$ then
\[
    \nabla\phi(R)^T\Delta \ge 0.
\]
Note that this means that $\nabla \phi(R)$ is a supporting hyperplane of $S$ at
$R$ if $\phi(R) = 1$. If the trader is trading some amount of asset $i$ for
asset $j$, \ie, $\Delta_i > 0$ and $\Delta_j < 0$ with all other entries zero,
we have
\[
    (\nabla\phi(R))_j\Delta_j \le -(\nabla\phi(R))_i\Delta_i,
\]
or, rewriting further,
\[
    \Delta_j \le \frac{(\nabla\phi(R))_i}{(\nabla\phi(R))_j}(-\Delta_i).
\]
Where equality can be achieved in the limit as the trade becomes small. We can
therefore interpret the quantity $(\nabla\phi(R))_i/(\nabla\phi(R))_j$ as the
price of token $j$ with respect to token $i$, and we can interpret the vector 
$\nabla \phi(R)$ as a vector of prices, up to a scaling factor determined by the
numeraire.

\paragraph{Discussion.} This shows that a number of results which hold `only'
for homogeneous trading functions, such as those
of~\cite{frongillo2023axiomatic,angerisNotePrivacyConstant2021}, are fully
general and hold for all CFMMs. Indeed, we do not need to assume homogeneity at
all as it may always be derived for a trading set satisfying some basic
conditions given above. Additionally, the direct connection to constant
function market makers comes from the fact that any trader may change the
reserves to some $R' \in \reals_+$ so long as
\[
    \phi(R') \ge \phi(R) = 1,
\]
where we assume that $\phi(R) = 1$ is a `starting condition' on the level set.
Of course, no trader would ever take $\phi(R') > 1$, since otherwise there
exists some dominating trade $\tilde R' \le R'$ with at least one inequality
holding strictly; \ie, the trader would tender less (or get more) of at least
one token and still have a feasible trade. So, in general, we have that, for
any `reasonable' action,
\begin{equation}\label{eq:cfmm}
    \phi(R') = \phi(R),
\end{equation}
where $R'$ is the new set of reserves, after a trade has been made, and $R$ is
the original set of reserves. Equation~\eqref{eq:cfmm} is the defining equation
for path-independent constant function market makers, explaining both their
name and the direct connection to the reachable set defined here.
(See~\cite{angeris2022constant} for more.)

\subsubsection{Uniqueness of canonical trading function}
We call this trading function \emph{canonical} since it is unique up to a
scaling constant. In fact, this function is unique if the function is scaled
such that the reachable set corresponds to its 1-superlevel set.

\paragraph{Proof.} To see this,
let $\phi$ and $\tilde \phi$ be two trading functions that are consistent and
yield the same reachable set $S$; \ie,
\[
    S = \{R \in \reals_+^n \mid \phi(R) \ge \alpha\} = \{R \in \reals_+^n \mid \tilde \phi(R) \ge \beta\},
\]
where $\alpha, \beta > 0$. (If $\alpha = 0$ then, since $\phi$ is homogeneous
and nondecreasing, we have that $\phi(R) \ge 0$, which would imply that its
reachable set is all of $\reals^n_+$, and similarly for $\tilde \phi$.) This
is the same as
\[
    \{R \mid \phi(R)/\alpha \ge 1\} = \{R \mid \tilde \phi(R)/\beta \ge 1\},
\]
so we will overload notation by writing $\phi$ for $\phi/\alpha$ and $\tilde
\phi$ for $\tilde\phi/\beta$, with the understanding that these differ by a
proportionality constant. Now, we will show that $\phi = \tilde \phi$.
To see this, start with the case that $R$ satisfies $\phi(R) > 0$ and $\tilde \phi(R) > 0$,
then
\[
    \phi\left(\frac{R}{\phi(R)}\right) = 1,
\]
so $R/\phi(R) \in S$ and we then have, by definition of $\tilde \phi$,
\[
    \frac{\tilde \phi(R)}{\phi(R)} = \tilde \phi\left(\frac{R}{\phi(R)}\right) \ge 1.
\]
Repeating the steps above with $\phi$ and $\tilde \phi$ swapped yields
\[
    \phi(R) = \tilde \phi(R),
\]
when $\phi(R) > 0$ and $\tilde \phi(R) > 0$. Now, if $\phi(R) = 0$ then 
\[
    \phi(tR) = t\phi(R) = 0,
\]
so $tR \not \in S$ for any $t > 0$. This means that
\[
    t \tilde \phi(R) = \tilde \phi(tR) < 1,
\]
again by definition of $\tilde \phi$, or, that $\tilde \phi(R) < 1/t$ for any
$t > 0$, so $\tilde \phi(R) = 0$. Repeating these steps where $\phi$ is swapped
with $\tilde \phi$ implies that $\phi(R) = 0$ only when $\tilde \phi(R) = 0$.
This gives the final result that $\phi = \tilde \phi$, or that the canonical
function is unique up to scaling constants.

\subsubsection{Examples}
In this subsection, we show the canonical trading function for Uniswap and
Uniswap v3. We also derive the canonical trading function for Curve~\cite{egorov2019stableswap}
in appendix~\ref{app:curve}.

\begin{figure}
    \centering
    \begin{tikzpicture}
    \begin{axis}[
        xlabel={$R_1$},
        ylabel={$R_2$},
        zlabel={$\lambda$},
        domain=0:2,
        y domain=0:2,
        zmin=0,
        zmax=2,
        mesh/cols=40,
        xtick=\empty,
        ytick=\empty,
    ]
    \addplot3[
        surf, 
        colormap={mycolormap}{rgb=(0.2,0.2,0.2) rgb=(1,1,1)},
        samples=40, 
        domain=0:2
    ] {sqrt(x*y)};
    \addplot3[black, ultra thick, samples y=0, domain=0.5:2] ({x},{1/x},{1})
    node[pos=0.5, above right, thick, black] {$\phi(R)$};
    \end{axis}
\end{tikzpicture}
    \hfill
    \begin{tikzpicture}
    \begin{axis}[
        xlabel={$R_1$},
        ylabel={$R_2$},
        xmin=0, xmax=5,
        ymin=0, ymax=5,
        legend pos=outer north east,
        xtick=\empty,
        ytick=\empty,
    ]
    \foreach \k in {0.1, 0.5, 1.5, 2} {
        \edef
        \temp{
            \noexpand
            \addplot[smooth, gray, thick, samples=500, domain=0.01:4.99,] ({x},{\k^2/x});
        }\temp
    }

    \addplot[smooth, very thick, samples=500, domain=0.01:4.99,] ({x},{1/x});

    \foreach \k in {2.5, 3, 3.5, 4, 4.5} {
        \edef
        \temp{
            \noexpand
            \addplot[smooth, gray, thick, samples=500, domain=1:4.99,] ({x},{\k^2/x});
        }\temp
    }
    
    \end{axis}
\end{tikzpicture}
    \caption{Left: the liquidity cone for Uniswap, with the level set defined by the
    trading function $\phi(R) = \sqrt{R_1 R_2} = 1$ shown. Right: each $\lambda$-level
    set of the surface looks like the boundary of the set of reachable reserves 
    (see figure~\ref{fig:uni-reserves}). The trading function $\phi$ is highlighted.}
    \label{fig:uni-liquidity-cone}
\end{figure}
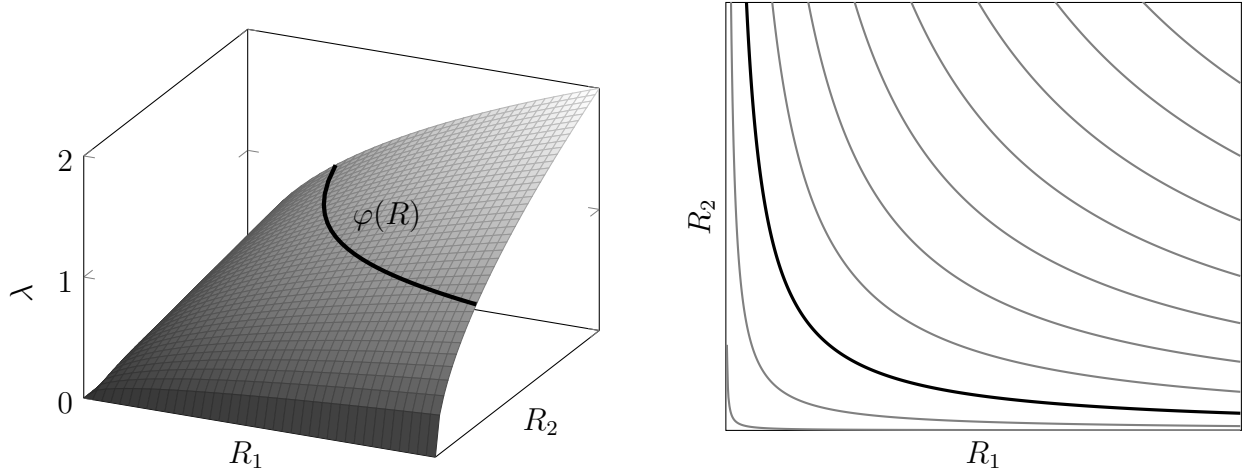

\paragraph{Uniswap.} Starting with the usual example of Uniswap, we have
that
\[
    S = \{R \in \reals_+^2 \mid R_1R_2 \ge k\}.
\]
The liquidity cone for Uniswap is given by
\begin{equation}\label{eq:uni-liquidity-cone}
    K = \{(R, \lambda) \in \reals^3 \mid R_1R_2 / \lambda^2 \ge k \in S, ~ \lambda > 0\}
\end{equation}
so, the canonical trading function~\eqref{eq:canonical-reach} can be written
\[
    \phi(R) = \sup \{\lambda > 0 \mid R_1R_2/\lambda^2 \ge k\},
\]
when $R \in \reals_+^2$ (and zero otherwise). This gives the canonical
trading function
\[
    \phi(R) = \sqrt{\frac{R_1R_2}{k}},
\]
which is evidently concave, nondecreasing, and 1-homogeneous, with
\[
    \phi(R) \ge 1 \quad \text{if, and only if} \quad R_1R_2 \ge k,
\]
as required. The liquidity cone and canonical trading function are shown in
figure~\ref{fig:uni-liquidity-cone}.

\paragraph{Uniswap v3.} We can also do the same for Uniswap v3, which has a
quasiconcave trading function given by
\[
    \psi(R) = (R_1 + \alpha)(R_2 + \beta).
\]
Since $\psi$ is strictly increasing in the positive orthant,
it suffices only to find the (positive) root of
\[
    (R_1/\lambda + \alpha)(R_2/\lambda + \beta) = k,
\]
which is a simple quadratic. The resulting canonical trading
function is unfortunately more complicated:
\begin{equation}\label{eq:uni-hom}
    \phi(R) = \frac12 \left(\frac{\beta R_1 + \alpha R_2 + \sqrt{(\beta R_1 + \alpha R_2)^2 + 4(k-\alpha\beta)R_1R_2}}{k - \alpha\beta}\right).
\end{equation}
This function is evidently homogeneous and strictly increasing since $k >
\alpha\beta$. Concavity is more difficult due to the square root term, but we
show it directly in appendix~\ref{app:concavity}. A good exercise is to show
that the canonical trading function $\phi$ in~\eqref{eq:uni-hom} has $\phi(R)
\ge 1$, if, and only if, $(R_1 + \alpha)(R_2 + \beta) \ge k$.

\subsection{Dual cone and portfolio value function}
In this section, we will look at an equivalent characterization of the
liquidity cone $K$, called the dual cone. The characterizations are equivalent
since the liquidity cone $K$ is convex. Indeed, we will show that this dual
cone has a very tight relationship with the portfolio value function, and leads
to a simple proof of the equivalence of (consistent) portfolio value functions
and (canonical) trading functions in that every portfolio value function has a
corresponding trading function, and vice versa, which was originally derived
in~\cite{angeris2023replicating}.

\subsubsection{Dual cone}
The \emph{dual cone} of a cone $K \subseteq \reals^{n+1}$ is defined as
\[
    K^* = \{(c, \eta) \in \reals^{n+1} \mid c^TR + \eta \lambda \ge 0, ~\text{for all}~ (R, \lambda) \in K\}.
\]
While this definition holds for any cone $K$, for the remainder of this
section, we will be working with the case that $K$ is the liquidity cone of a
CFMM with reachable set $S$, as defined the previous subsection. 

\paragraph{Intuition.} In a very general sense, the dual cone $K^*$ is simply
another (dual) representation of the original liquidity cone, $K$, in that the
dual of $K^*$, defined as $(K^*)^* = K$, as $K$ is closed and convex. (For more
information on conic duality, we refer the reader to
appendix~\ref{app:conic-duality}.) We will use this fact to give a simple proof
that the trading function and the portfolio value function (to be introduced
later in this section) are two views of the same underlying object. 

\paragraph{Basic properties.} First, note that $K^*$ is always a closed, convex
cone as it can be written as the intersection of closed hyperplanes, and, by
definition, we have $0 \in K^*$. Additionally, we have that
\begin{equation}\label{eq:nonneg-dual}
    K^* \subseteq \reals_+^n \times \reals
\end{equation}
since $K \supseteq \reals_{++}^n \times \{0\}$, from the previous section. (To
see this, use the definition of $K^*$.) Finally, we may write the dual cone in
terms of only the reachable set $S$. We have that $(c, \eta) \in K^*$ if, and
only if,
\[
    c^TR + \eta \lambda \ge 0, ~ \text{for all} ~ (R, \lambda) \in K,
\]
by definition. But this latter statement is true if, and only if it is
true for all $\lambda > 0$, since $K$ is the closure over 
the set defined in~\eqref{eq:liquidity-cone}; \ie, $(c, \eta) \in K^*$ if, and only if,
\[
    c^TR + \eta \lambda \ge 0, ~ \text{for all} ~ (R, \lambda) \in K, ~ \lambda > 0.
\]
Rearranging the inequality gives that $c^T(R/\lambda) + \eta \ge 0$, and note that,
by definition of $K$, we have that $(R, \lambda) \in K$ with $\lambda > 0$ only when
$R/\lambda \in S$. This means that $(c, \eta) \in K^*$ if, and only if,
\begin{equation}\label{eq:alt-dual}
    c^T\tilde R + \eta \ge 0, ~ \text{for all} ~ \tilde R \in S.
\end{equation}
This particular rewriting of $K^*$ will be useful in what follows.


\subsubsection{Portfolio value function}\label{sec:pv}
Much in the same way that we defined
the trading function, we may define the \emph{portfolio value function} as
\begin{equation}\label{eq:pv-def}
    V(c) = \sup \{-\eta \mid (c, \eta) \in K^*\}.
\end{equation}
This function has the following interpretation: given an external
market with prices $c \in \reals_+^n$ (\ie, anyone may trade asset $i$ for asset $j$
at a fixed price $c_i/c_j$) then $V(c)$ corresponds to the total
value of reserves after arbitrage has been performed. In particular,
$V(c)$ is the optimal value of the problem,
\begin{equation}\label{eq:pv-problem}
    \begin{aligned}
        & \text{minimize} && c^TR\\
        & \text{subject to} && R \in S,
    \end{aligned}
\end{equation}
with variable $R \in \reals^n$, where $S$ is the reachable set.

To see this, note that $(c, \eta) \in K^*$ if, and only if,
\[
    c^TR \ge -\eta, ~ \text{for all} ~ R \in S,
\]
from the previous characterization of $K^*$ given in~\eqref{eq:alt-dual}.
The claim follows by applying the definition of $V$ in~\eqref{eq:pv-def}.

\paragraph{Properties.} 
From the optimization problem formulation~\eqref{eq:pv-problem}, we see that
$V$ is clearly nonnegative and nondecreasing since for any $R \in S$, we have
that $R \ge 0$. The function $V$ is also concave because it is the partial
maximization of the concave function
\[
    f(c, \eta) = -\eta - I(c, \eta),
\]
over $\eta$, where $I$ is the indicator function of the convex set $K^*$,
defined as $I(c, \eta) = 0$ if $(c, \eta) \in K^*$ and $+\infty$,
otherwise. Finally, we see that $V$ is homogeneous since for $\alpha > 0$,
$V(\alpha c)$ is the optimal value of the problem
\[
    \begin{aligned}
        & \text{minimize} && \alpha c^TR\\
        & \text{subject to} && R \in S.
    \end{aligned}
\]
Since $\alpha$ is a constant, this value is clearly $\alpha V(c)$.

\paragraph{Consistency.} We say a portfolio value function is \emph{consistent}
if it is concave, homogeneous, and nondecreasing, which we know is true for any
function $V$ derived from a reachable set $S$. Of course, every consistent
portfolio value function also defines a dual cone:
\[
    K^* = \{(c, \eta) \mid V(c) + \eta \ge 0\},
\]
which can be easily verified to be a convex cone that is contained in $K^*
\subseteq \reals_+^n \times \reals$ using the fact that $V$ is consistent,
so we may convert from portfolio value functions to dual cones directly.

\paragraph{Examples.} Using~\eqref{eq:uni-liquidity-cone}, we can write the
dual cone for Uniswap 
\[
    K^* = \{
        (c, \eta) \mid c^TR + \eta \lambda \ge 0, 
        ~ \text{for all} ~ R_1 R_2 \ge k \lambda^2, ~ \lambda > 0 
    \}.
\]
We can simplify this expression via a few observations. First, we must have $c
\ge 0$, from~\eqref{eq:nonneg-dual}. Second, because $c \ge 0$, if $\eta \ge 0$
then $(c, \eta)$ is clearly in $K^*$. The interesting case is then when $c \ge
0$ but $\eta < 0$. In this case, we must have that
\[
    c^T R \ge (-\eta) \sqrt{\frac{R_1R_2}{k}},
\] 
since $\lambda$ can take any value between $0$ and $\sqrt{R_1 R_2 / k}$.
Rearranging, we have that
\[
    c_1 x + c_2 x^{-1} \ge (-\eta) / \sqrt{k},
\]
where $x = \sqrt{R_1 / R_2}$. Minimizing the left hand side over $x > 0$ means
that this inequality is true if, and only if,
\[
    2\sqrt{c_1c_2} \ge -\eta/\sqrt{k},
\]
so the dual cone for Uniswap is
\[
    K^* = \{(c, \eta) \in \reals_+^2 \times \reals \mid 2\sqrt{kc_1c_2} + \eta \ge 0\}.
\]
The portfolio value function can almost be read off from the definition:
\begin{equation}\label{eq:univ2-pv}
    V(c) = 2\sqrt{kc_1c_2},
\end{equation}
which is evidently concave, homogeneous, and nondecreasing.

As a more complicated example, we'll derive the portfolio value function for a
Uniswap v3 `tick'. In this case, it's easier to work directly from the
optimization problem~\eqref{eq:pv-problem}. For convenience, let $c = (p, 1)$
and note that we can recover the general case using the homogeneity of $V$, as
$V(\tilde c) = \tilde c_2 V(\tilde c_1/\tilde c_2, 1)$. Then,
\[
    V(p, 1) = \inf_{R \ge 0}\left\{pR_1 + R_2 \mid (R_1 + \alpha)(R_2 + \beta) \ge k\right\}.
\]
Any profit maximizing trader will ensure that the inequality holds with
equality (\ie, the solution is at the boundary of the set $S$). After
substitution, we have a simple convex function that is minimized either at a
point $R > 0$ with gradient zero or at the boundary. We can conclude that
\[
    V(p, 1) = \begin{cases}
        pk/\alpha - \beta               & p < \beta^2/k \\
        2\sqrt{pk} - (\alpha + \beta)   & \beta^2/k \le p \le k/\alpha^2 \\
        pk/\beta - \alpha               & k/\alpha^2 < p.
    \end{cases}
\]
Note the similarity of this expression to the previous~\eqref{eq:univ2-pv} when
the price is within a particular range. This range corresponds exactly to the
`tick' interval in Uniswap v3~\cite{adams2021uniswap}.

\subsubsection{Replicating market makers}\label{sec:rmm}
In this subsection we show how to convert directly between the portfolio value
function and a canonical trading function (and vice versa). This shows that,
indeed, every canonical trading function has an equivalent consistent portfolio
value function, and, in a sense, each of these functions is a different `view'
of the same underlying object.

\paragraph{Trading function to portfolio value.}
Assuming $\phi$ is a canonical trading function, as defined
in~\eqref{eq:canonical-cone}, then we may write the portfolio value function as
\[
    V(c) = \inf_{R > 0}\left(\frac{c^TR}{\phi(R)}\right).
\]
To see this, note that the definition of $K^*$ is that $(c, \eta) \in K^*$
when
\[
    c^TR + \eta \lambda \ge 0, ~\text{for all} ~ (R, \lambda) \in K.
\]
Minimizing the left hand side over $\lambda > 0$ gives that $(c, \eta) \in K^*$
only when
\[
    c^TR + \eta \phi(R) \ge 0, ~ \text{for all} ~ R \in \reals_+^n,
\]
by definition of $\phi(R)$. Using a basic limiting argument, we may replace
$R \in \reals_+^n$ with $R \in \reals_{++}^n$, which implies that $\phi(R) > 0$
by positive reachability~\eqref{eq:positive-reachability}, so we have that
\[
    -\eta \le \frac{c^TR}{\phi(R)}, ~ \text{for all} ~ R \in \reals_{++}^n,
\]
or, equivalently, that $(c, \eta) \in K^*$ if, and only if,
\[
    -\eta \le \inf_{R > 0}\left(\frac{c^TR}{\phi(R)}\right).
\]
Applying the definition of $V(c)$, given in~\eqref{eq:pv-def}, gives the final
result.

\paragraph{Trading function from portfolio value.} It is also possible
to show that we can recover a canonical trading function from a given
portfolio value function. To see this, note that
\begin{equation}\label{eq:pv-to-trade}
    \phi(R) = \inf_{c > 0}\left(\frac{c^TR}{V(c)}\right),
\end{equation}
is a concave (as it is the minimization of a family of affine functions over $R$), 
homogeneous, and nondecreasing trading function. We can easily show
that, if $K^*$ corresponds to the dual of a liquidity cone $K$, and $V$ is the
corresponding portfolio value function, then $\phi(R)$ corresponds to its
canonical trading function.

From a nearly-identical argument to the previous, replacing the definition of
$\phi$ with that of $V$, we have that $(\tilde R, \tilde \lambda) \in (K^*)^*$
if, and only if,
\[
    \tilde \lambda \le \inf_{c > 0}\left(\frac{c^T\tilde R}{V(c)}\right).
\]
Since $K$ is a liquidity cone (by assumption) it is therefore closed and
convex, so we have that $(K^*)^* = K$; \cf, appendix~\ref{app:conic-duality}.
Finally, maximizing over $\tilde \lambda$ and using the definition of $\phi$
given in~\eqref{eq:canonical-cone}:
\[
    \phi(\tilde R) = \inf_{c > 0}\left(\frac{c^T\tilde R}{V(c)}\right),
\]
where $\phi$ is the canonical trading function for $K$.

\paragraph{Example.} To complete the cycle, we convert the portfolio value
function of Uniswap back to its canonical trading function. From above,
\[
\begin{aligned}
    \phi(R) = \inf_{c > 0}\left(\frac{c^TR}{V(c)}\right)
    &= \frac{1}{2\sqrt{k}}\inf_{c > 0}\left(\sqrt{\frac{c_1}{c_2}}R_1 + \sqrt{\frac{c_2}{c_1}}R_2\right) \\
    &= \frac{1}{2\sqrt{k}}\inf_{x > 0}\left(xR_1 + x^{-1}R_2\right) \\
    &= \sqrt{\frac{R_1R_2}{k}},
\end{aligned}
\]
where we recognized $xR_1 + x^{-1}R_2$ as a convex function and minimized by
simply applying the first order optimality conditions.

\paragraph{Interpretation.} There is a nice interpretation for
equation~\eqref{eq:pv-to-trade} which is that the quotient $c^TR/V(c)$ denotes
the \emph{leverage} or the `lambda' of the portfolio $R \in \reals_+^n$ at
price $c$, where $V(c)$ denotes the true value of the CFMM holdings at this
price. We may then view the trading function $\phi(R)$ as the lowest possible
leverage over all possible prices. The inequality $\phi(R) \ge 1$, which
defines the reachable set, says that the leverage must be at least 1 in order
for the reserves to lie in the set.

\paragraph{Connection to RMMs.} There is a connection to the original
result of~\cite{angeris2023replicating} by noting that the trading function 
presented there is defined, using the portfolio value function $V$, as
\[
    \phi^0(R) = \inf_{c > 0}\left(c^TR - V(c)\right) = -I_{(K^*)^*}(R, 1),
\]
where $I_{(K^*)^*}$ is the indicator function for the dual cone of the dual
cone, $(K^*)^*$. Since $(K^*)^* = K$ then $\phi^0(R) \ge 0$ if, and only if,
$(R, 1) \in K$, which happens if, and only if, $R \in S$, as required.

\paragraph{Discussion.} From the above, we have that every consistent portfolio
value leads to a canonical trading function. This method gives a general
procedure for going from one to the other. Additionally, since we know
$((K^*)^*)^* = K^*$, then we know that, starting from any consistent portfolio
value function $V$, converting it to a trading function $\phi$, and then
converting back results in the same $V$ we started with, which shows that the
mapping is indeed invertible.

\subsubsection{Composition rules}
We will denote $S_V \subseteq \reals_+^n$ as the reachable set corresponding to
the portfolio value function $V$. (We will see how to construct this explicitly
in what follows.)

\paragraph{Composition rules for portfolio value.} Given consistent portfolio
value functions, there are a number of possible ways these could be `combined'.
The first is by scaling: if $V$ is consistent, then certainly $\alpha V$ is
consistent. If both $V$ and $V'$ are consistent, then $V+V'$ is consistent,
and, finally if $A$ is a nonnegative orthogonal matrix, and $V$ is consistent,
then $V\circ A^T$ is consistent. We will show that these operations correspond
to natural operations over the reachable sets corresponding to the portfolio
value functions. 

\paragraph{Recovering the reachable set.} We may recover the reachable set from a
given portfolio value function since we know that, for given portfolio value
function, its liquidity cone, which we will denote $K_V \subseteq
\reals^{n+1}$ is
\[
    K_V = \{ (R, \lambda) \in \reals^{n+1} \mid c^TR -\lambda V(c) \ge 0 ~ \text{for all} ~ c \ge 0\}.
\]
Clearly, this cone is convex, closed, and satisfies $K_V \subseteq
\reals_+^{n+1}$ from the fact that $V$ is consistent. Additionally,
since we may define a reachable set from a liquidity cone as $S_V = \{R \mid
(R, 1) \in K_V\}$, then this is the same as saying
\begin{equation}\label{eq:reachable-pv}
    S_V = \{R \mid c^TR \ge V(c) ~ \text{for all} ~ c \ge 0\}.
\end{equation}
It remains to be verified that $S_V$ is a valid reachable set, but this follows
from the properties of $K_V$ outlined above. (Another way to see this is to
note that $c^TR\ge V(c)$ if, and only if, $c^TR/V(c)\ge 1$ for all $c > 0$,
\ie, $\phi(R) \ge 1$ using~\eqref{eq:pv-to-trade}. Since $\phi$ is a canonical
trading function, then $S_V$ is a reachable set.)

\paragraph{Scaling.} It is not hard to see that
\[
    S_{\alpha V} = \alpha S_V,
\]
for any $\alpha \ge 0$ by using the definition of $S_V$ and the fact that
$V$ is homogeneous.

\paragraph{Addition.} Similarly, addition over the portfolio value functions
`commutes' over the reachable sets; \ie,
\[
    S_{V + V'} = S_V + S_{V'}.
\]
The direction $S_{V + V'} \subseteq S_V + S_{V'}$ is easy to show by definition.
On the other hand, since $S_{V+V'}$ is a closed convex set, if $R \not \in S_{V+V'}$
then there exists a strictly separating hyperplane $c \in \reals_+^n$ with 
\[
    c^TR < c^T\tilde R, ~ \text{for all} ~ \tilde R \in S_{V+V'}.
\]
Taking the infimum of the right hand side and using~\eqref{eq:pv-problem} gives
\[
    c^TR < (V+V')(c) \le c^T\tilde R + c^T \tilde R' ~ \text{for all} ~ \tilde R \in S_V, ~ \tilde R' \in S_{V'},
\]
which means that $R \not \in S_V + S_{V'}$. Here, the last inequality follows
from the fact that $V(c) \le c^T\tilde R$ for all $\tilde R \in S_V$ and
similarly for $S_{V'}$. Putting both statements together gives that $S_{V+V'} =
S_V + S_{V'}$.

\paragraph{Nonnegative projection.} Given some nonnegative matrix $A \in
\reals_+^{m \times n}$ that is also an orthogonal matrix, \ie, $A^TA = I$, then
\[
    S_{V\circ A^T} = AS_V +\reals_+^m,
\]
where $(V\circ A^T)(c) = V(A^Tc)$. This follows nearly immediately from the
definition of $A$ and~\eqref{eq:reachable-pv}.

\paragraph{Intersection.} There is a natural question then, as to what the
intersection of reachable sets corresponds to. Clearly, we have
\[
    S_V \cap S_{V'} = \{R \mid c^TR \ge V(c) ~ \text{and} ~ c^TR \ge V'(c)\}.
\]
Of course, this implies that
\[
    S_V \cap S_{V'} = S_{\max\{V, V'\}},
\]
where the $\max$ is taken pointwise. Note that this does not correspond
to a natural operation on the portfolio value functions as the pointwise
maximum of two concave functions is not necessarily concave. (Take, for
example, $V(p, 1) = \sqrt{p}$ and $V'(p, 1) = p$.) Let, on the other hand,
$\tilde V$ be the (pointwise) smallest concave function with $\tilde V \ge V$
and $\tilde V \ge V'$, then indeed we have
\[
    S_V \cap S_{V'} = S_{\tilde V},
\]
and it is not hard to show that $\tilde V$ is consistent.

\paragraph{Discussion.} This also gives another proof of the composition rules
presented for the reachable sets since we may always recover a consistent
portfolio value from any reachable set. In this sense, we may think of the
portfolio value function and the reachable set as two objects with a `natural
homomorphism' under which scaling, addition, nonnegative projection, and
intersection are all preserved.

\subsection{Connection to prediction markets} Prediction markets are a type of
market which attempts to elicit the beliefs of players about the probability
that certain events occur. These markets have a rich academic history, stemming
back since at least the 50s~\cite{mccarthy1956measures} and, until the
relatively recent paper~\cite{frongillo2023axiomatic}, a connection between
such markets and CFMMs was not known, except in some very special cases.  This
section restates and simplifies the results of~\cite{frongillo2023axiomatic} in
this framework. We differ in the notion of `histories' for path-independent
CFMMs which is implicitly included in this framework and discussed later in
this paper, as a general result in~\S \ref{sec:path-independence}.

\paragraph{Cost functions.} A \emph{cost function} is defined as a 
function $C: \reals^n \to \reals \cup \{+\infty\}$ such that
\begin{enumerate}
    \item The function $C$ is convex, nondecreasing
    \item It is translation invariant, $C(q + \alpha \ones) = C(q) + \alpha$
    \item It is somewhere finite; \ie, there is at least one $q$ for which $C(q)$ is finite
\end{enumerate}
Note that we do not require the function to be differentiable, only
subdifferentiable in its domain. This means that there might be many
probabilities which are consistent with the market's predictions, but includes
differentiability as the special case where there is only one.

\paragraph{Example.} One particular, important example is the logarithmic
market scoring rule, or LMSR, which has cost function
\[
    C(q)= b\log\left(\sum_{i=1}^n \exp\left(\frac{q_i}{b}\right)\right),
\]
where $b > 0$ is some given constant. This function is clearly nondecreasing
and finite. It is also convex as it is a log-sum-exp~\cite[\S3.1.5]{cvxbook} function with
affine precomposition. It is not hard to see that this function is also
translation invariant using the definition, which means that this function is,
indeed, a reasonable cost function.

\paragraph{Mechanics.} The mechanics of a prediction market are that any player
may buy any of $n$ possible mutually exclusive outcomes. Every player will be
paid a dollar for each share of outcome $i$ they hold if outcome $i$ occurs at
some future time. All other outcomes will have a value of zero. A player who
wishes to buy $\delta \in \reals^n$ shares and must pay a cost of $C(q +
\delta) - C(q)$, where $q$ is the current set of outstanding shares that have
been sold to all players. (Negative values of $\delta_i$ means that the player
is selling back $\delta_i$ shares to the market.) The outstanding shares are
then updated to $q \gets q + \delta$.

\paragraph{Interpretation.} Let's say the prediction market begins with some
outstanding shares $q_0$, and a player has beliefs $p \in \reals_+^n$ about the
probability of each event such that $p_i$ corresponds to the probability of the
$i$th event occurring. The player can then maximize her expected profit (under
her distribution of beliefs) by solving
\[
    \begin{aligned}
        & \text{maximize} && p^Tq - (C(q_0 + q) - C(q_0)),
    \end{aligned}
\]
with variable $q$. We note that the optimal value of this problem, call it
$E(p)$ for expected payoff at probabilities $p$, is tightly related to the Fenchel
conjugate of $C$, since
\[
    E(p) = C^*(p) + C(q_0) - p^Tq_0.
\]
The optimality conditions for this problem are that
\[
    p \in \partial C(q_0 + q^\star),
\]
where $q^\star$ is the solution to the optimization problem. This means that,
if the market has some outstanding shares given by $q_0$ then we may interpret
$\partial C(q_0)$ as the set of probabilities consistent with the market's
belief about the event. 

\paragraph{CFMM to cost function.} Given a reachable set $S$, we can construct
a cost function:
\begin{equation}\label{eq:reach-to-cost}
    C(q) = \min\{\alpha \ge 0\mid \alpha\ones - q \in S\}.
\end{equation}
(This was first observed by~\cite{frongillo2023axiomatic}.) We may also define the
cost function in terms of the liquidity cone as
\[
    C(q) = \max\{\beta \ge 0 \mid (\ones - \beta q, \beta) \in K\}.
\]
This function is a cost function since it is evidently translation invariant by
definition, and is nondecreasing since $S$ is upward closed. The function is
finite at $0$ since $S$ is nonempty: if $R \in S$, then $0 \le C(0) \le \max_i
R_i$. Additionally, this function is convex as it is the partial minimization
of the convex function
\[
    f(\alpha, q) = \alpha + I(\alpha, q),
\]
over $\alpha \in \reals$, where $I(\alpha, q) = 0$ if $\alpha \ones - q \in S$
and $\alpha \ge 0$, and is $+\infty$ otherwise. (This set indicator is convex
as it is the indicator function for the intersection of convex sets.) 

\paragraph{Example.} Recall that Uniswap has the trading set 
\[
    S = \{ R \in \reals_+^2 \mid R_1 R_2 \geq k\}
\]
Using~\eqref{eq:reach-to-cost}, we have the cost function
\[
    C(q) = \min \{\alpha \geq 0 \mid (\alpha - q_1)(\alpha - q_2) \geq k\}.
\]
The cost function is the positive root of the quadratic over $\alpha$, the same
as was found in~\cite{frongillo2023axiomatic}:
\[
    C(q) = \frac{q_1 + q_2}{2} + \frac{1}{2}\sqrt{(q_1 - q_2)^2 + k}.
\]
We can easily verify that this function is finite, translation invariant, and
convex (by noting that the square root term can be expressed as the $\ell_2$
norm of the vector $(q_1-q_2, \sqrt{k})$). The fact that it is nondecreasing
can be seen by showing that its gradient is everywhere nonnegative.

\paragraph{Cost function to CFMM.} Any cost function $C$ defines a CFMM
by defining its reachable set as,
\begin{equation}\label{eq:cost-to-reach}
    S = \{R \in \reals_+^n \mid C(-R) \le 0\}.
\end{equation}
This $S$ is indeed a reachable set as (a) the function $C$ is nondecreasing by
assumption, so $S$ is upward closed, (b) it is convex as $C$ is convex, and (c)
it is nonempty since $C(q)$ is finite for some $q \in \reals^n$, so $C(q - C(q)
\ones ) = 0$ by translation invariance, and therefore $C(q)\ones - q \in S$.
We may write its canonical trading function using~\eqref{eq:canonical-reach}:
\[
    \phi(R) = \sup\{\lambda > 0 \mid C(-R/\lambda) \le 0\}.
\]

\paragraph{Equivalence.} If the cost function $C$ is constructed from a CFMM
with reachable set $S$, as in~\eqref{eq:reach-to-cost}, then it is not hard to
show that~\eqref{eq:cost-to-reach} yields exactly this set $S$. To see this,
note that, by definition~\eqref{eq:reach-to-cost}, we have that $C(-R) \le 0$,
if, and only if, $\alpha\ones + R \in S$ for all $\alpha \ge 0$; letting
$\alpha = 0$ gives that $R \in S$. On the other hand, if $R \in S$ then, $R +
\alpha \ones \in S$ for every $\alpha \ge 0$, by upward closedness, so $C(-R)
\le 0$ and the sets are equivalent.

\paragraph{Example.} 
The logarithmic market scoring rule (LMSR) has the cost function
\[
    C(q)= b\log\left(\sum_{i=1}^n \exp\left(\frac{q_i}{b}\right)\right).
\]
We may define its trading set, using~\eqref{eq:cost-to-reach}, as
\[
    S = \left\{R \in \reals_+^n ~ \middle|~ \sum_{i=1}^n \exp(-R_i/b) \le 1\right\}.
\]
The corresponding canonical trading function is
\[
    \phi(R) = \sup\left\{\lambda > 0 ~\middle|~  \sum_{i=1}^{n} \exp(-R_i/\lambda b) \le 1 \right\}   
\]
This function has no closed form solution but can be solved numerically as a
univariate root-finding problem. Since $C$ is strictly increasing, the positive
root is unique and can be found efficiently using the methods discussed
in~\S\ref{sec:trading-function}.

\subsection{Liquidity provision}\label{sec:liquidity-provision}
As in~\cite{angeris2022constant}, we discuss liquidity provision in the case
where the trading function $\phi$ is homogeneous. This is, of course, fully
general as we may assume that $\phi$ is a consistent trading function.
 
\paragraph{Liquidity providers.} An agent, called a \emph{liquidity provider}
can add or remove assets from the CFMM's reserves $R$. When an agent adds
liquidity, she adds a basket $\Psi \in \reals^n_+$ to the reserves, resulting
in the updated reserves $R^+ = R + \Psi$. When an agent removes liquidity, she
removes a basket $\Psi \in \reals^n_+$ with $\Psi \le R$ from the reserves,
resulting in the updated reserves $R^+ = R - \Psi$. When adding (or removing)
to reserves, the agent receives (tenders) an IOU. This IOU gives the agent a
pro-rata share of the reserves based on the amount of value the agent added and
the total amount of value in the pool. We describe the exact mechanism for
liquidity addition (and removal) below.

\paragraph{Liquidity change condition.}
The main condition for adding and removing liquidity is that the asset prices
must not change after the removal, or addition, of liquidity. More
specifically, we must have that the prices, as given
in~\S\ref{sec:marg-prices}, at the new reserves, $R^+ \in \reals_+^n$,
must be equal, up to a scalar multiple, to those at the original reserves,
$R \in \reals_+^n$. Written out, this gives the following condition:
\[
    \nabla \phi(R^+) = \alpha\nabla \phi(R),
\]
where $\phi$ is the canonical trading function and $\alpha > 0$ is some
positive constant. Since $\phi$ is homogeneous, we have that, for any
$\alpha > 0$, $\nabla \phi(\alpha R) = \alpha\nabla \phi(R)$. We conclude that
$\Psi = \nu R$ for $\nu > 0$ is a valid liquidity change for any $\nu > 0$
(where we must have $\nu \le 1$ for liquidity removal). 
Note that scaling $R$ to $\alpha R$ corresponds
exactly to scaling the reachable set by $\alpha$.

\paragraph{Liquidity provider share weights.} The CFMM additionally maintains a
table of all liquidity providers and their corresponding share weights,
representing the fraction of the reserves that each liquidity provider owns. We
denote these weights as $w \in \reals_+^N$, where $N$ is the number of
liquidity providers, and enforce that they sum to one, \ie, $\sum_{i=1}^N w_i =
1$. These weights are updated whenever a liquidity provider adds or removes
liquidity, or when the number of liquidity providers $N$ changes.

\paragraph{Value of the reserves.} Let $V = p^TR$ be the value of the CFMM
reserves at price $p$. After adding liquidity $\nu R$, the value of the
reserves is now
\[
    V^+ = p^TR^+ = (1 + \nu)p^TR = (1 + \nu)V.
\]
For removing liquidity, we replace $\nu$ with $-\nu$. The fractional change in
reserve value is
\[
    (V^+ - V)/V^+ = \nu/(1 + \nu).
\]

\paragraph{Liquidity provider share update.} When liquidity provider $j$ adds
or removes liquidity, all the share weights are adjusted pro-rata based on the
change of value of the reserves, which is the value of the basket she adds or
removes. The fractional change in reserve value is $\nu/(1 + \nu)$. Thus, after
adding liquidity, the change in share weights is
\[
    w_i^+ = \begin{cases}
        w_i /(1 + \nu) + \nu / (1 + \nu) & i = j \\
        w_i / (1 + \nu) & i \neq j.
    \end{cases}
\]
For removing liquidity, we replace $\nu$ with $-\nu$ and add the constraint that
$\nu \le w_j$.

\paragraph{Portfolio value.} We note that, since liquidity providers own a
pro-rata share with weight $w_i$ of the total pool value, we may view each
liquidity provider's position as `independent'. In particular, there is no
distinction between many liquidity providers pooling their assets together into
a single CFMM versus every liquidity provider having their own CFMM instance
and owning all of the assets of their particular instance. (There may be
practical differences, however, owing to the fact that users may prefer to
trade with a subset of these for a variety of reasons, such as gas fees.)

\section{Single trade}
We consider in this section the general CFMM case, which potentially includes
fees and is therefore not necessarily path-independent. (We show the connection
to the previous fee-free case later in the section.)

\subsection{Trading set}
Much in the same way as the previously-defined reachable set, we will define
the \emph{trading set} $T \subseteq \reals^n$, which is any set $T$ satisfying
the following properties:
\begin{enumerate}
    \item The set $T$ is closed and convex
    \item The zero trade is included, $0 \in T$
    \item The set $T$ is downward closed; \ie, if $\Delta \in T$ and $\Delta' \le \Delta$ then $\Delta' \in T$
\end{enumerate}
An additional requirement that will be useful in the composition
rules presented later, but is not strictly required
for most of the statements below is that: there exists $R \in \reals_+^n$
such that
\begin{equation}\label{eq:finite-tender}
    R - T = \{R - \Delta \mid \Delta \in T\}\subseteq \reals^n_+.
\end{equation}
This corresponds to the statement that the CFMM can only tender a finite amount
of some asset (in `usual' CFMMs, this would be the available reserves) which is
upper bounded by the quantity $R \ge 0$. One could imagine a mechanism that is
allowed to mint unbounded amounts of tokens may
violate~\eqref{eq:finite-tender}.

\paragraph{Set up.} In this set up, we have a trader who wishes to trade with
the CFMM. This trader can suggest to trade a basket of tokens $\Delta \in
\reals^n$, where positive values denote what the trader receives from the CFMM
and negative values denote what the trader tenders to the CFMM. The CFMM
accepts this trade (and tenders or receives the stated amounts) only when
$\Delta \in T$. (If $\Delta \not \in T$, the trader receives and tenders
nothing to the CFMM.) The state of the CFMM is then updated based on the
accepted trade $\Delta$ (and a rejected trade does not change the state), but,
for now, we will only consider the single-trade case and elide discussion of
this state update until later. We assume only that the trading set at the
current state is known and accessible to the trader.

\paragraph{Interpretation.} The conditions imposed on the trading set all have
relatively simple interpretations. The convexity of the trading set means that,
as users trade more of a token, they receive marginally less (or, at least, no
more) than they otherwise would by trading less. The fact that the zero trade
is included means that a user is allowed to not trade. Finally, the
downward-closedness of the set means that the CFMM will accept more of a given
token, all else being equal; \ie, a trader is allowed to `overpay' for a given
trade, and this new trade is still valid. The final optional condition can be
interpreted as: a CFMM has a finite amount of assets that it is allowed to
tender. While not strictly a requirement, we will need it for a technical
condition we present later.

\begin{figure}
    \centering
    \begin{tikzpicture}[scale=1.0]
    \begin{axis}[
        axis on top=true,
        xmin=-4, xmax=2,
        ymin=-4, ymax=2,
        axis lines=center,
        xlabel={$\Delta_1$},
        ylabel={$\Delta_2$},
        grid=major,
        legend pos=outer north east,
        samples=200,
        domain=-5:5,
        xtick=\empty,
        ytick=\empty
        ]                    
        \addplot[black, very thick, domain=0:1.99] {2*x/(x - 2)};
        \addplot[black, very thick, domain=0:1.99] ({2*x/(x - 2)}, x);
        \addplot[gray!90, domain=0:1.99, fill, opacity=0.5, draw=none] ({2*x/(x - 2)}, x) -| (-5,-5) -- cycle;
        \addplot[gray!90, domain=0:1.99, fill, opacity=0.5, draw=none] {2*x/(x - 2)} -| (-5,-5) -- cycle;

        \addplot[black, very thick, domain=0:0.99] {x/(x - 1)};
        \addplot[black, very thick, domain=0:0.99] ({x/(x - 1)}, x);
        \addplot[gray!30, domain=0:0.99, fill, opacity=1, draw=none] ({x/(x - 1)}, x) -| (-5,-5) -- cycle;
        \addplot[gray!30, domain=0:0.99, fill, opacity=1, draw=none] {x/(x - 1)} -| (-5,-5) -- cycle;
        
        \node at (-2.5, -2.5) {$T$};
    \end{axis}
\end{tikzpicture}
    \hfill
    \begin{tikzpicture}[scale=1.0]
    \begin{axis}[
        axis on top=true,
        xmin=-4, xmax=2,
        ymin=-4, ymax=2,
        axis lines=center,
        xlabel={$\Delta_1$},
        ylabel={$\Delta_2$},
        grid=major,
        legend pos=outer north east,
        samples=200,
        domain=-4:4,
        xtick=\empty,
        ytick=\empty
        ]
        \addplot[black, very thick, domain=0:0.5] {x/(x - 1)};
        \addplot[black, very thick, domain=0:0.5] ({x/(x - 1)}, x);
        \addplot[gray!30, domain=0:0.5, fill, opacity=0.6, draw=none] ({x/(x - 1)}, x) -| (-4,-4) -- cycle;
        \addplot[gray!30, domain=0:0.5, fill, opacity=0.6, draw=none] {x/(x - 1)} -| (0.5, -1) -- (0.5, -4) -- (-4,-4) -- cycle;
        \node at (-2.5, -2.5) {$T$};

        \draw[black, very thick] (-1, 0.5) -- (-4, 0.5);
        \draw[black, very thick] (0.5, -1) -- (0.5, -4);
    \end{axis}
\end{tikzpicture}
    \caption{Left: the trading set for Uniswap (without fees) for $R = (1, 1)$
    (light gray) and $R = (2,2)$ (dark gray). Right: the trading set for Uniswap 
    v3.}
    \label{fig:uni-trading-set}
\end{figure}
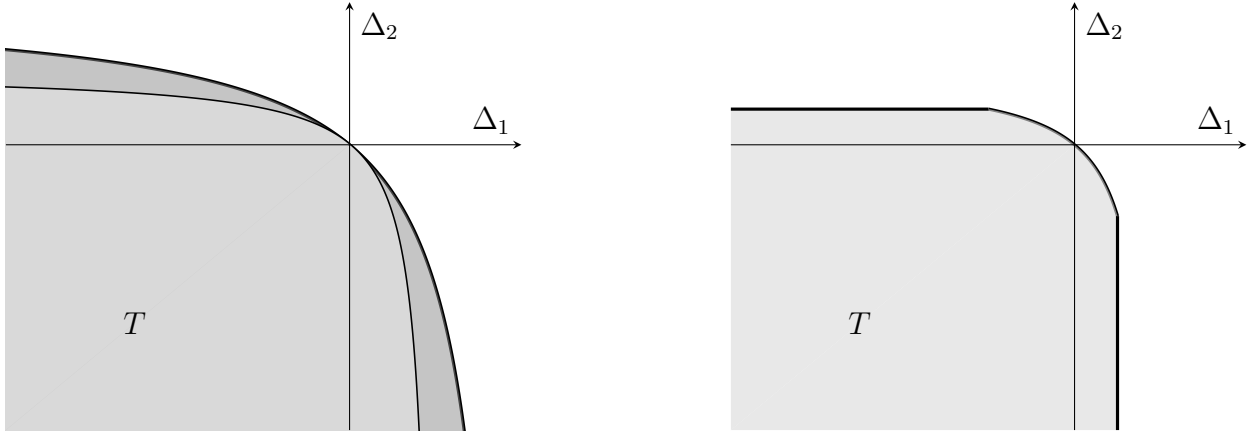

\paragraph{Example.} A basic example of a trading set is Uniswap
with fees. In this particular case, the current state
of the CFMM is given by some reserves $R \in \reals_+^2$,
and the trading set is
\[
    T = \{\Delta^+ - \Delta^- \mid (R_1+\gamma\Delta^-_1 - \Delta^+_1)_+(R_2 + \gamma \Delta^-_2 - \Delta^+_2)_+ \ge R_1R_2, ~\Delta^-, \Delta^+ \in \reals_+^n\},
\]
where $0 \le \gamma \le 1$ is the fee parameter and is set by the CFMM at
creation time. We show this set in figure~\ref{fig:uni-trading-set}. This 
writing is bit verbose and difficult to parse, but the construction is very 
similar to the original, given in the fee-free case above. Because of this, it 
is often nicer to work directly with a functional form, which we describe below.

\paragraph{Quasiconcavity and fees.} Given any nondecreasing, quasiconcave
function, with nonnegative domain $\psi:\reals_+^n \to \reals$ (much like the
previous) we can define a trading set with fee $0 \le \gamma \le 1$ and
reserves $R \in \reals_+^n$,
\begin{equation}\label{eq:quasiconcave-trading}
    T = \{\Delta^+ - \Delta^- \mid \psi(R+\gamma \Delta^- - \Delta^+) \ge \psi(R), ~ \Delta^-, \Delta^+ \in \reals_+^n \}.
\end{equation}
(We will, by convention, let the function $\psi$ take on $-\infty$ for values
outside of its domain.) Clearly $0\in T$ and $T$ is downward closed as $\psi$
is nondecreasing. The set is evidently convex as it is an affine transform of
the set
\[
    \{(\Delta^+, \Delta^-) \in \reals_+^n\times \reals_+^n \mid \psi(R+\gamma \Delta^- - \Delta^+) \ge \psi(R)\},
\]
which is a convex set by the quasiconcavity of $\psi$. Closedness is trickier,
but follows from the fact that the valid choices of $\Delta^+$ are
in some compact set, $0 \le \Delta^+ \le R$.

\paragraph{Composition rules.}\label{sec:composition-trade} The composition rules are nearly
identical in both statement and proof to those of the reachable set.
Given trading sets $T, T' \subseteq \reals^n_+$
\begin{enumerate}
    \item Trading sets may be added; \ie, $T + T'$ yields a trading set
    \item Trading sets may be scaled, so $\alpha T$ is a trading set for any $\alpha > 0$
    \item Taking the intersection $T \cap T'$ preserves the trading set property
    \item Applying a nonnegative linear transformation $A \in \reals_+^{k\times n}$
        and adding all dominated trades,
        \[
            AT- \reals_+^k,
        \]
        preserves the trading set property
\end{enumerate}
These composition rules similarly lead to the notion of an aggregate CFMM
mentioned previously, in the single-trade case, which is especially useful in
the case of Uni v3 as we will show later. The only technical condition appears
when adding trading sets: to ensure that the resulting trading set is closed,
it suffices to ensure that all CFMMs can tender only a finite amount of assets,
as in condition~\eqref{eq:finite-tender}.

\subsection{Trading cone and dual}
Much in the same way as we have defined the liquidity cone, we define
the \emph{trading cone} as
\[
    K = \cl \{(\Delta, \lambda) \in \reals^{n+1} \mid \Delta/\lambda \in T, ~ \lambda > 0\}.
\]
This cone plays a similar role to the liquidity cone, except in the single
trade case. Indeed, many of the constructions we have shown previously for the
liquidity cone will apply in a similar form to the trading cone.
(We have overloaded notation as we will make no further reference to the
liquidity cone.)

\paragraph{Trading function.} From a nearly identical argument to the previous
we may define a homogeneous, nondecreasing, but \emph{convex} (instead of
concave) trading function
\[
    \phi(\Delta) = \min\{\lambda \ge 0 \mid (\Delta,\lambda) \in K\},
\]
or, equivalently,
\begin{equation}\label{eq:gen-trading-function}
    \phi(\Delta) = \inf\{\lambda > 0 \mid \Delta/\lambda \in T\}.
\end{equation}
such that
\[
    T = \{\Delta \in \reals^n \mid \phi(\Delta) \le 1\}.
\]
(Here, we define the $\min$ and $\inf$ of an empty set to be $+\infty$ for
convenience.) The difference in sign from the definition in the path
independent case in~\S\ref{sec:trading-function} comes from the fact that the set
$T$ is downward (rather than upward) closed, since we are taking the
perspective of the trader, rather than the CFMM or its liquidity providers. In
this case, the function $\phi$ is similarly canonical and rational traders will
always tender trades $\Delta$ such that
\[
    \phi(\Delta) = 1,
\]
hence, again, the name `constant function market maker.'

\paragraph{Example.} Perhaps the simplest example of this type of function is,
unsurprisingly, Uniswap. Using the quasiconcave function
definition~\eqref{eq:quasiconcave-trading} of the trading set, we have,
for $R \in \reals_+^2$,
\[
    \psi(R_1, R_2) = R_1R_2
\]
with some fee $0 \le \gamma \le 1$. For a given proposed trade, $\Delta\in
\reals^n$, we can decompose $\Delta$ into its positive and negative
parts $ \Delta = \Delta^+ - \Delta^-$ with $\Delta^-, \Delta^+ \ge 0$ and
disjoint support $\Delta^-_i\Delta^+_i = 0$ for each $i=1, 2$. Using the
definition of the trading function, we look for the smallest $\lambda \ge 0$
such that
\[
    \left(R_1 + \gamma\frac{\Delta^-_1}{\lambda} - \frac{\Delta^+_1}{\lambda}\right)\left(R_2 + \gamma\frac{\Delta^-_2}{\lambda} - \frac{\Delta^+_2}{\lambda}\right) \ge R_1R_2.
\]
With some basic rearrangements, we find
\[
    \phi(\Delta) = \frac{(\Delta^+_1 - \gamma\Delta^-_1)(\Delta^+_2 - \gamma\Delta^-_2)}{R_1(\gamma\Delta^-_2 - \Delta^+_2) + R_2(\gamma\Delta^-_1 - \Delta^+_1)}.
\]
This trading function is homogeneous since the numerator is a homogeneous
quadratic while the denominator is homogeneous. We can similarly verify that,
as expected from the previous discussion, it is convex and nondecreasing by
writing it in the following form:
\[
    \phi(\Delta) = \frac{1}{-(R_1/(\gamma\Delta^-_1 - \Delta^+_1) + R_2/(\gamma\Delta^-_2 - \Delta^+_2))}.
\]
Since the denominator is nonnegative, concave, and nonincreasing (in
$\Delta^+$), then $\phi$ must be nonnegative, convex, and nondecreasing (in
$\Delta^+$). Since $\Delta = \Delta^+ - \Delta^-$, directly verifying fact
that $\phi$ is nondecreasing and convex in $\Delta$ requires one more
step, which we leave to the reader as a useful exercise. (Of course, we know
that both of these already follow from the construction of the trading function
in~\eqref{eq:gen-trading-function} and the fact that $\psi$ is quasiconcave and
nondecreasing, as is the case for all such trading functions.)

\paragraph{Bounded liquidity.} in a similar way to the previous section, we
know that, since $0 \in T$, then
\[
    (-\reals_+^n, 0) \subseteq K.
\]
We say the trading set has \emph{bounded liquidity in asset $i$} if
the supremum
\[
    \sup\{\Delta_i \mid \Delta \in T\} = \Delta^\star_i,
\]
is achieved at some $\Delta^\star \in T$. This has the interpretation that
there is a finite basket of assets such that we receive all possible amount of
asset $i$ from the CFMM. We say a trading set has \emph{bounded liquidity} if
it has bounded liquidity for each asset $i=1, \dots, n$. Examples of bounded
liquidity CFMMs include Uniswap v3 (see figure~\ref{fig:uni-trading-set}) and those with linear trading functions.
These bounded liquidity CFMMs are useful since arbitrage can be easily computed
in many important practical cases; see~\cite[\S3]{diamandis2023efficient}
for more.

\paragraph{Arbitrage cone.} In a similar way to the previous, we will define
the dual cone for the trading cone $K\subseteq\reals^{n+1}$ as
\[
    K^* = \{(c, \eta) \mid c^T\Delta + \eta \lambda \ge 0, ~ \text{for all} ~ (\Delta, \lambda) \in K\}.
\]
By downward closedness and the fact that $0 \in T$, it is not hard to show that
$K^*\subseteq (-\reals_+)^n\times \reals_+$. Minimizing over the left hand side
of the inequality gives another definition, based on the trading function:
\[
    K^* = \{(c, \eta) \mid c^T\Delta + \eta\phi(\Delta) \ge 0, ~ \text{for all} ~ \Delta \in \reals^n\}.
\]
Some care has to be taken when interpreting this expression if $\phi(\Delta) =
\infty$ when $\eta = 0$, based on the original definition of $K$, but this is
an informative exercise for the reader. 

\paragraph{Relation to arbitrage.} Much like the portfolio value function, we
write the \emph{arbitrage function}, $\arb: \reals^n \to \reals$, for the
trading set $T$ as
\begin{equation}\label{eq:arb-def}
    \arb(c) = \sup\{c^T\Delta \mid \Delta \in T\}.
\end{equation}
Note that if $c_i < 0$ for any $i$ then $\arb(c) = +\infty$ by the
downward-closedness of $T$, so we may generally assume that $c \ge 0$. This
function has the following interpretation: if there is an external market with
prices $c \in \reals_+^n$, this is the maximum profit that an arbitrageur could
derive by trading between the external market and the CFMM. This function is
convex (as it is the supremum of a family of functions that are affine in $c$),
nondecreasing over $c \ge 0$, and homogeneous. We may write this function in
terms of the dual cone as
\[
    \arb(c) = \inf\{\eta \mid (-c, \eta) \in K^*\},
\]
by a nearly-identical argument to that of the portfolio value function
in~\S\ref{sec:pv}. This function will be very useful in the routing problem
that follows. Additionally, from a very similar argument to~\S\ref{sec:rmm},
the arbitrage function and the trading function are equivalent representations
in that we may derive one from the other by setting
\[
    \phi(\Delta) = \sup_{\arb(c) > 0} \left(\frac{c^T\Delta}{\arb(c)}\right),
\]
and
\[
    \arb(c) = \sup_{\phi(\Delta) > 0} \left(\frac{c^T\Delta}{\phi(\Delta)}\right).
\]
From before, note the suprema in this equation, versus the infima in the
previous. For examples of such arbitrage functions for some common constant
function market makers, see~\cite[app.\ A]{diamandis2023efficient}.

\begin{figure}
    \centering
    \begin{tikzpicture}[scale=1]
    \begin{axis}[
      xmin=-5, xmax=5,
      ymin=-5, ymax=5,
      axis lines=center,
      xlabel={$\Delta_1$},
      ylabel={$\Delta_2$},
      grid=major,
      legend pos=outer north east,
      samples=200,
      domain=-5:5,
      xtick=\empty,
      ytick=\empty
      ]
  
      \node at (-2.5, -2.5) {$T$};
  
      \addplot[black, very thick, domain=0:0.99] {2*x/(x - 1)};
  
      \addplot[black, very thick, domain=0:0.99] ({2*x/(x - 1)}, x);
  
      \addplot[gray!30, domain=0:0.99, fill, opacity=0.3, draw=none] ({2*x/(x - 1)}, x) -| (-5,-5) -- cycle;
      \addplot[gray!30, domain=0:0.99, fill, opacity=0.3, draw=none] {2*x/(x - 1)} -| (-5,-5) -- cycle;
  
      \addplot[black, very thick, domain=0:2.4] {2*x};
      \addplot[black, very thick, domain=0:4.8] {0.5*x};
      \addplot[gray!100, fill, opacity=0.3, draw=none] (0,0) -- (4.8, 2.4) -- (2.4, 4.8) -- cycle;
      \node at (2.4, 2.4) {$C(0)$};
  
      \addplot[black, dashed, domain=-2:2] {-2*x};
      \addplot[black, dashed, domain=-4:4] {-x/2};

    \end{axis}
  \end{tikzpicture}
    \caption{The trading set for Uniswap with fees (notice that the set is 
    kinked at $0$) and the corresponding no-trade cone.}
    \label{fig:univ2-marginal-price-cone}
\end{figure}
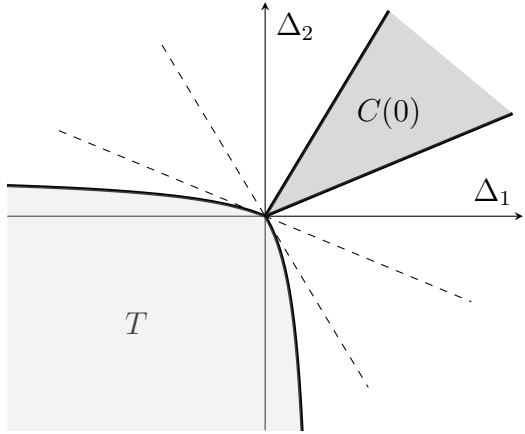

\paragraph{Marginal prices.} We can view the supporting hyperplanes of $T$
at some $\Delta$ as the set of \emph{marginal prices} at trade $\Delta$.
We write this set as
\begin{equation}\label{eq:marg-prices}
    C(\Delta) = \bigcap_{\Delta' \in T} \{\nu \in \reals^n \mid \nu^T(\Delta' - \Delta) \le 0\}.
\end{equation}
Note that this set is a closed convex cone as it is the intersection of closed
convex cones and is always nonempty as $0 \in C(\Delta)$. We can write the cone
$C(\Delta)$ using the trading function as
\begin{equation}\label{eq:marg-subgrad}
    C(\Delta) = \bigcup_{\lambda \ge 0} \lambda \partial(\phi(\Delta)),
\end{equation}
whenever $\phi(\Delta) = 1$ and the subdifferential is defined. As we will soon
see, the cone $C(0)$ will be called the \emph{no-trade cone}. This is a
generalization of the \emph{no-trade interval}~\cite{diamandis2023efficient} in
the case where $n \ge 2$. We show this cone for Uniswap in figure~\ref{fig:univ2-marginal-price-cone}.

The proof of the equivalence~\eqref{eq:marg-subgrad} can be shown in two steps,
one for the forward inclusion, and one for the reverse. From the statement, we
have $\Delta$ with $\phi(\Delta) = 1$. Now let $g \in \partial \phi(\Delta)$,
then
\[
    \phi(\Delta)+g^T(\Delta' - \Delta) \le \phi(\Delta'),
\]
by definition of the subgradient $g$. Letting $\Delta' \in T$ means that
$\phi(\Delta') \le 1$ by definition and $\phi(\Delta) = 1$ by the previous,
so
\[
    g^T(\Delta' - \Delta) \le \phi(\Delta') - 1 \le 0,
\]
for every $\Delta' \in T$, which means that $g \in C(\Delta)$.
Multiplying both sides of this inequality by $\lambda \ge 0$
then means that $\lambda g \in C(\Delta)$, or that
\[
    \bigcup_{\lambda \ge 0} \lambda \partial(\phi(\Delta)) \subseteq C(\Delta).
\]

For the other direction, using the definition of $C$ in~\eqref{eq:marg-prices},
we can see that
$g \in C(\Delta)$ if, and only if, $\Delta$ is a maximizer of the following
optimization problem:
\[
    \begin{aligned}
        & \text{maximize} && g^T\Delta'\\
        & \text{subject to} && \phi(\Delta') \le 1,
    \end{aligned}
\]
with variable $\Delta' \in \reals^n$. Using the optimality conditions of this
problem, we know that $\Delta$ is a maximizer if, and only if, there exists
some $\lambda \ge 0$ such that
\[
    0 \in -g + \lambda \partial\phi(\Delta),
\]
or, equivalently, if, and only if, $g \in \lambda \partial\phi(\Delta)$ for some
$\lambda \ge 0$, which shows the reverse inclusion. It also shows that
\[
    C(\Delta) \subseteq \reals_+^n,
\]
since $\phi$ is nondecreasing, so it subgradients must be nonnegative. The
optimality conditions are necessary and sufficient by Slater's condition~\cite[\S5.2.3]{cvxbook}, 
since $\phi(-\ones) = 0 < 1$ and
$-\ones$ is in the interior of the domain of $\phi$.

\paragraph{Marginal price composition.} Given $\Delta_i \in T_i$ for $i=1,
\dots, m$, we have that
\begin{equation}\label{eq:price-cone}
    C\left(\sum_{i=1}^m \Delta_i\right) = \bigcap_{i=1}^m C_i(\Delta_i),
\end{equation}
where $C_i$ is the cone of marginal prices for CFMM $i$ while $C$ is the cone
of marginal prices for the aggregate CFMM
\[
    \tilde T = \sum_{i=1}^m T_i.
\]
This is easy to see from the definitions of $C$ and $\tilde T$.

\paragraph{Connection to arbitrage.} Note that $\Delta$ is a solution to the
arbitrage problem at price $c$, \ie, $c^T\Delta = \arb(c)$ if, and only if,
\[
    c \in C(\Delta),
\]
which follows by using the definition of $\arb$ and $C$. In other words, the
arbitrage problem is solved at any trade which changes the prices to match
those of the external market with prices $c \in \reals^n_+$. We say there
is \emph{no arbitrage at price $c$} if the zero trade is a solution, \ie,
\[
    c \in C(0).
\]
Equivalently, we may view this as the case where the CFMM's prices are
consistent with those of the external market. Alternatively, there
is a direct connection between the marginal prices, arbitrage, and
the dual cone $K^*$:
\[
    c \in C(\Delta), ~\text{if, and only if,} ~ (-c, c^T\Delta) \in K^*.
\]
We can see this since $c \in C(\Delta)$ if, and only if, for all $\Delta' \in T$, we have
\[
    c^T\Delta' \le c^T\Delta.
\]
But $(\Delta', \lambda') \in K$ with $\lambda' > 0$ if and only if $\Delta'/\lambda' \in T$
so
\[
    \frac{c^T\Delta'}{\lambda'} \le c^T\Delta.
\]
Multiplying both sides by $\lambda' > 0$ and using a limiting argument shows that this
is true, if, and only if, for all $(\Delta', \lambda') \in K$ we have
\[
    c^T\Delta' \le \lambda' c^T\Delta,
\]
which is the same as saying $(-c, c^T\Delta) \in K^*$.

\subsection{Routing problem}
The routing problem takes a number of possible CFMMs $i=1, \dots, m$, each 
trading a subset of $n_i$ tokens out of the universe of $n$ tokens, and seeks to 
find the best possible set of trades, \ie, those maximizing a given utility 
function $U: \reals^n \to \reals \cup \{-\infty\}$. We assume that $U$ is
concave and increasing (\ie, we assume all assets have value with potentially
diminishing marginal returns). We use infinite values of $U$ to encode 
constraints; a trade $\Psi$ such that $U(\Psi) = -\infty$ is unacceptable to the
trader. See~\cite[\S5.2]{angeris2022optimal} for examples, including liquidating or
purchasing a basket of tokens and finding arbitrage.

We denote the trade we make with the $i$th CFMM by $\Delta_i$ and this CFMM's trading 
cone by $K_i \subseteq \reals^{n_i+1}$.
We also introduce matrices $A_i \in \reals^{n\times n_i}$ which map the `local' 
basket of $n_i$ tokens for CFMM $i$ to the global universe of $n$ tokens.
This construction is similar to the construction of aggregate CFMMs in~\S\ref{sec:composition-rules},
but here we focus on the trade vectors and not the trading sets. The net trade
is simply
\[
    \Psi = \sum_{i=1}^m A_i \Delta_i.
\]
The optimal routing problem is then the problem of finding a set of valid trades
with each market that maximize the trader's utility:
\[
    \begin{aligned}
        & \text{maximize} && U(\Psi)\\
        & \text{subject to} && \Psi = \sum_{i=1}^m A_i\Delta_i\\
        &&& (\Delta_i, 1) \in K_i, \quad i=1, \dots, m.
    \end{aligned}
\]
The variables here are the net trade $\Psi \in \reals^n$ and the trades
$\Delta_i \in \reals^{n_i}$. Note that, by definition of the trading cone
$K_i$, we have that $\Delta_i \in T_i$ if, and only if, $(\Delta_i, 1) \in
K_i$.

\paragraph{Other interpretations.} If $A_i = I$, \ie, if all CFMMs trade the
same tokens, then this problem is equivalent to
\[
    \begin{aligned}
        & \text{maximize} && U(\tilde \Delta)\\
        & \text{subject to} && \tilde \Delta \in \tilde T,
    \end{aligned}
\]
where $\tilde T = \sum_{i=1}^m T_i$, which is another trading set, by the
composition rules given in~\S\ref{sec:composition-trade}. While this rewriting
seems silly, it tells us that we may consider routing through a network of
CFMMs as trading with one `large' CFMM.  The optimality conditions for this
problem are that
\[
    0 \in \partial (-U)(\tilde \Delta^\star) + \tilde C(\Delta^\star),
\]
and $\tilde \Delta^\star \in \tilde T$. From~\eqref{eq:price-cone} we know that
$\tilde C$ is the intersection of each individual price cone, so, using the
definition of $\tilde T$, we get
\[
    0 \in \partial(-U)\left(\sum_{i=1}^m A_i\Delta_i^\star\right) + \bigcap_{i=1}^m C_i(\Delta_i^\star),
\]
and $\Delta_i^\star \in T_i$, which are exactly the optimality conditions
we would get from considering the original routing problem. The case where
$A_i$ are general nonnegative orthogonal matrices is slightly more involved,
but is ultimately very similar.

\paragraph{Dual problem.} From conic duality (\cf, appendix~\ref{app:conic-duality}),
we know that the dual problem can be written as
\[
    \begin{aligned}
        & \text{minimize} && \bar U(\nu) + \ones^T\eta\\
        & \text{subject to} && (-A_i^T\nu, \eta_i) \in K_i^*, \quad i=1, \dots, m,
    \end{aligned}
\]
where the variables are $\nu \in \reals^n$ and $\eta \in \reals^m$.
Partially minimizing over each $\eta_i$ and using the definition of the optimal
arbitrage function, we have that this problem is equivalent to
\[
    \begin{aligned}
        & \text{minimize} && \bar U(\nu) + \sum_{i=1}^m \arb_i(A_i^T\nu),
    \end{aligned}
\]
where $\arb_i$ is the optimal arbitrage function for the $i$th trading set.
This is exactly the dual problem used in the decomposition method of~\cite{diamandis2023efficient}.
This problem has a beautiful interpretation: the optimal trades are exactly
those which result in a price vector $\nu$ that minimizes the total arbitrage profits
that the user would receive if we interpret $\bar U(\nu)$ as the maximum utility
that could be received by trading with an external market with price $\nu$.

\subsection{Path independence}\label{sec:path-independence}
In this subsection, we show the connection between the path-independent CFMMs,
presented in the previous section, and the `general' CFMMs presented in this
one. 

\paragraph{Mechanics of trading.} In a CFMM, as stated previously, we have some
state, which is given by the reserves $R \in \reals_+^n$. The current trading
set, defined as $T(R) \subseteq \reals^n$ has the same properties given at the
beginning of this section. (We implicitly included the relationship between the
trading set and the reserves in the previous section as the reserves could be
considered fixed for a single trade.) The CFMM then accepts or rejects any
proposed trade $\Delta \in \reals^n$ based on whether $\Delta \in T(R)$. If
this is the case, then the CFMM accepts the trade, updating its reserves to $R
\to R - \Delta$ (as it pays out $\Delta_i$ to the trader from its reserves for
$\Delta_i > 0$ and vice versa) and making the new trading set $T(R - \Delta)$.
If the trade is rejected then the reserves are not updated and the trading
set remains as-is.

\paragraph{Sequential feasibility.} From before, we say a trade $\Delta$ is
feasible if $\Delta \in T(R)$. We say a sequence of trades $\Delta_i \in
\reals^n$ for $i=1, \dots, m$ is \emph{(sequentially) feasible} if
\[
    \Delta_{i} \in T(R - (\Delta_1 + \dots + \Delta_{i-1})).
\]
for each $i=1, \dots, m$.

\paragraph{Reachability.} We say some reserves $R'$ are \emph{reachable}
from some initial set of reserves $R$ if there is a sequence of feasible
trades $\Delta_i \in \reals^n$, for $i=1, \dots, m$, such that
\[
    R' = R - (\Delta_1 + \dots + \Delta_m).
\]
In other words, $R'$ is reachable from $R$ if there is a sequence of feasible
trades that takes us from reserves $R$ to reserves $R'$

\paragraph{Path independence.} We say a CFMM is \emph{path independent} if, for
any reserves $R$ and for any trade $\Delta$ satisfying $\Delta \in T(R)$, we
have
\begin{equation}\label{eq:path-independence}
    \Delta' \in T(R-\Delta) \quad \text{if, and only if}, \quad \Delta + \Delta' \in T(R).
\end{equation}
In English: a CFMM is path independent if there is no difference between
performing trades sequentially versus in aggregate, if the trades are
sequentially feasible. (We may apply induction to this definition to get the
more `general-seeming' case that applies to any finite sequence of feasible
trades.)

\begin{figure}
    \centering
    \begin{tikzpicture}[scale=0.9]
    \begin{axis}[
        axis on top=true,
        xmin=-4, xmax=1,
        ymin=-4, ymax=1,
        axis lines=center,
        xlabel={$\Delta_1$},
        ylabel={$\Delta_2$},
        grid=major,
        legend pos=outer north east,
        samples=200,
        domain=-4:4,
        xtick=\empty,
        ytick=\empty
        ]
        \addplot[black, very thick, domain=0:0.99] {x/(x - 1)};
        \addplot[black, very thick, domain=0:0.99] ({x/(x - 1)}, x);
        \addplot[gray!30, domain=0:0.99, fill, opacity=0.5, draw=none] ({x/(x - 1)}, x) -| (-4,-4) -- cycle;
        \addplot[gray!30, domain=0:0.99, fill, opacity=0.5, draw=none] {x/(x - 1)} -| (-4,-4) -- cycle;
        \node at (-2.5, -2.5) {$T(R')$};
    \end{axis}
\end{tikzpicture}
    \hfill
    \begin{tikzpicture}
    \begin{axis}[
        axis on top=true,
        xmin=-1, xmax=4,
        ymin=-1, ymax=4,
        axis lines=center,
        x label style={anchor=north},
        y label style={anchor=east},
        xlabel={$R_1$},
        ylabel={$R_2$},
        grid=major,
        legend pos=outer north east,
        samples=200,
        domain=-4:4,
        xtick=\empty,
        ytick=\empty
        ]

        \addplot[black, very thick, domain=0.01:3.99] {1/x};
        \addplot[gray!30, domain=0.01:4.99, fill, opacity=0.5, draw=none] (1/x, x) -| (4,4) -- cycle;
        \node at (2.5, 2.5) {$S = R' - T(R') $};
        \node[below left] at (1.0, 1.0) {$R'$};
        \fill (1.0, 1.0) circle[radius=2pt];
    \end{axis}
    \end{tikzpicture}
    \caption{The trading set $T(R')$ for Uniswap (left) and the corresponding
    reachable set $S$ (right).}
    \label{fig:uni-trading-and-reserve-sets}
\end{figure}
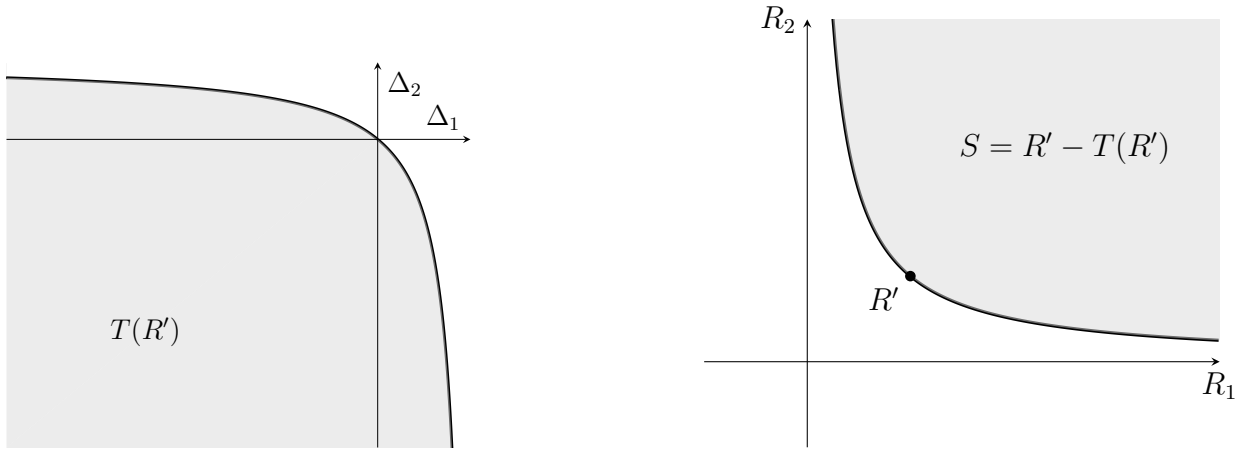

\paragraph{Reachable set.} If the CFMM is path independent, there exists a
fixed set $S\subseteq \reals^n$ (which, as we will soon see, corresponds
exactly to the reachable set of~\S\ref{sec:reachable-set}) such that every
trading set $T(R')$ can be written as
\begin{equation}\label{eq:feasible-set}
    T(R') = R' - S,
\end{equation}
for any reachable $R'$, starting from some reserves $R$. (Here, $R' - S = \{R' -
\tilde R \mid \tilde R \in S\}$.) Figure~\ref{fig:uni-trading-and-reserve-sets}
illustrates these sets for Uniswap.

\paragraph{Proof.} We will show that, whenever the CFMM is path independent,
then we will have that, for any reachable $R'$, $R - T(R) = R' - T(R')$. Setting
$S = R - T(R)$ will then suffice to satisfy~\eqref{eq:feasible-set}. Note
that it suffices to consider only $R'$ which are reachable in 1 step, since the
result follows by induction. That is, we will consider the case where $R' = R - \Delta$
for $\Delta \in T(R)$ and the general case follows by induction.

We can rewrite the path independence condition~\eqref{eq:path-independence} as
\[
    \Delta' \in T(R - \Delta) \quad \text{if, and only if}, \quad \Delta' \in T(R) - \Delta,
\]
or, equivalently,
\[
    T(R - \Delta) = T(R) - \Delta.
\]
The proof is then nearly obvious after this:
\[
    R' - T(R') = (R - \Delta) - T(R - \Delta) = R - \Delta - T(R) + \Delta = R - T(R).
\]
We may then set $S = R' - T(R') = R - T(R)$, such that~\eqref{eq:feasible-set} is satisfied
for any $R'$ reachable from $R$.

\paragraph{Conditions.} Note that the conditions on $T(R)$ will imply some conditions
on $S$. Indeed, since $T(R)$ is a closed convex set, then $S$ must also be. Similarly, since $0 \in T(R)$
then $R \in S$, so $S$ is nonempty, and, since we must have $R - \Delta \ge 0$ for any $\Delta \in T(R)$
then $S = R - T(R) \subseteq \reals_+^n$. Finally, since $T(R)$ is downward closed, then $S$ must be
upward closed, hence $S$ must be a reachable set as defined in~\S\ref{sec:reachable-set}.

\paragraph{Equivalence.} This is, of course, a bijection. We know that any path
independent CFMM with trading set $T(R)$ may be written as
\[
    T(R) = R - S,
\]
so long as $R \in S$. Additionally, if $S$ is a reachable set, then we must
have that $0 \in T(R)$, that $T(R)$ is closed and convex, and that $T(R)$ is
downward closed, making it a reasonable trading set. It is also
bounded~\eqref{eq:finite-tender} since $S \subseteq \reals_+^n$ so $R - T(R)
\subseteq \reals_+^n$.

\paragraph{Discussion.} In general, it is tempting to deal with histories of
trades among other objects when discussing CFMMs as these are dynamic systems
with some internal state that changes as trades are performed. The above shows
that, in the special case that the CFMM is path independent, we only need to
consider the reachable set, as this contains all properties needed to
completely describe the object in question. Indeed, the proof above shows that
a CFMM is path-independent if, and only if, it is completely described by a
reachable set meeting the conditions outlined in~\S\ref{sec:reachable-set}.

\section{Conclusion}
In this paper, we have shown that a general geometric perspective on constant
function market makers is relatively fruitful. Indeed, assuming only a small
number of `intuitive' conditions on sets, we have derived a number of
results---some known already, some not---which follow from almost purely
geometrical considerations. In some cases, we show that assumptions made in the
literature actually are unnecessary, and indeed are consequences of a subset of
the assumptions made. Examples of these include the homogeneity
of~\cite{angerisNotePrivacyConstant2021}, \cite{schlegel2023axioms},
and~\cite{frongillo2023axiomatic}. In others, we derive a new form of a known
result such as~\cite{angeris2023replicating} or~\cite{frongillo2023axiomatic}.
We suspect that there are a number of useful `geometric' interpretations to
other known results in the literature, but leave these for future work.

\bibliographystyle{alpha}
\bibliography{citations.bib}

\appendix
\section{A primer on conic duality}\label{app:conic-duality}
This appendix is intended as a (very short) primer on conic duality. We assume
basic familiarity with convex sets and the separating hyperplane theorem. For
far more, see~\cite[\S2.6]{cvxbook}.

\paragraph{Cones.} A \emph{cone} is a set $K \subseteq \reals^n$ such that, if $x
\in K$ then, for any $\alpha \ge 0$, we have $\alpha x \in K$. A \emph{convex
cone} is, as one would expect, a cone that is convex. More generally, convex
cones are closed under nonnegative scalar multiplication, \ie, if $x, y \in K$
and $\alpha, \beta \ge 0$, then
\[
    \alpha x + \beta y \in K.
\]
Basic examples of convex cones include the nonnegative real elements $\reals_+^n$
and the norm cones, given by
\[
    K_{\|\cdot\|} = \{(x, t) \in \reals^{n+1} \mid \|x\| \le t\}.
\]

\paragraph{Properties.} If $K$ is closed and nonempty then $0 \in K$. The
intersection of convex cones is a convex cone and scaling a convex cone results
in a convex cone. Finally, the (Cartesian) product of two convex cones is again
a convex cone.

\paragraph{Dual cone.} The \emph{dual cone} $K^*$ of a cone $K$ is defined
as
\[
    K^* = \{y \in \reals^n \mid y^Tx \ge 0, ~ \text{for all} ~ x \in K\}.
\]
In other words, the dual cone of $K$ is the set of all vectors which have nonnegative
inner product with every element in $K$. Since we can write
\[
    K^* = \bigcap_{x \in K} \{y \in \reals^n \mid y^Tx \ge 0\},
\]
then we can see that $K^*$ is a closed convex cone (even when $K$ is not).
For example, the dual cone of $\reals_+^n$ is $\reals_+^n$, while the norm
cone is
\[
    K_{\|\cdot\|}^* = \{(y, r) \in \reals^{n+1} \mid \|y\|_* \le r\},
\]
where $\|y\|_*$ is the dual norm of $y$, defined
\[
    \|y\|_* = \sup\{y^Tx \mid \|x\| \le 1\}.
\]

\paragraph{Properties.} By definition $0 \in K$, and, since $K^*$ can be
written as an intersection over $K$, then if $K' \subseteq K$, we have
\[
    K^* \subseteq K'^*.
\]
Additionally note that
\[
    (K+K')^* = K^*\cap K'^*,
\]
and
\[
    (K\times K')^* = K^*\times K'^*,
\]
all of which are simple exercises and follow from the definition of
the dual cone above.

\paragraph{Duality.} In a certain sense, we may view the dual cone $K^*$ as a
collection of \emph{certificates} that an element is not in $K$. More
specifically, if we have any $x \in \reals^n$ and we are given some $y \in K^*$
such that $y^Tx < 0$, then we are guaranteed that, indeed $x \not \in K$, by
definition of $K^*$. Conic duality gives the following guarantee for a nonempty, closed
convex cone $K$: for any $x \in \reals^n$, either $x \in K$, or there exists $y
\in K^*$ with $y^Tx < 0$, but not both. In other words, either some given
point $x$ it either belongs in the cone, or we can furnish a certificate,
using the dual cone, that it does not.

The reverse implication follows from the previous argument. (Note that this
implication requires no assumptions on $K$.) The forward implication will make
use of the convexity of $K$ to furnish a certificate. To see this, let $x \not
\in K$, then, since $K$ is convex closed, there exists a strict separating
hyperplane with slope $y \in \reals^n$ such that
\[
    x^Ty < z^Ty, ~ \text{for all} ~ z \in K.
\]
Since, for any $t \ge 0$ and $z \in K$, we have $tz \in K$, we therefore
know that, for any $z \in K$,
\[
    x^Ty/t < z^Ty,
\]
and sending $t \to \infty$ we then know 
\[
    z^Ty \ge 0 ~\text{for any} ~ z \in K,
\]
so $y \in K^*$. Finally, since $K$ is closed, then $0 \in K$ so
\[
    x^Ty < 0,
\]
completing the proof.

\paragraph{Dual of the dual.} Because of the previous, we now have the following
result: $K$ is exactly the set of vectors $x \in \reals^n$ such that $x$
has nonnegative inner product with every element of $K^*$; \ie, for which
we cannot furnish a certificate that $x \not \in K$. But, the set of vectors
which have nonnegative inner product with every element of $K^*$ is exactly
the dual cone of $K^*$, written $(K^*)^* = K^{**}$. This gives the following beautiful
relation for a nonempty, closed, and convex cone $K$:
\[
    K^{**} = K.
\]

\paragraph{Conic duality in optimization.} Most convex optimization problems can be
cast as conic optimization problems. The general form of such a problem is,
for some convex objective function $f: \reals^n \to \reals \cup \{\infty\}$
\[
    \begin{aligned}
        & \text{minimize} && f(x)\\
        & \text{subject to} && Ax = b\\
        &&& x \in K,
    \end{aligned}
\]
where the variable is $x \in \reals^n$, and the problem data are the closed
nonempty convex cone $K\subseteq\reals^n$, the matrix $A \in \reals^{m\times n}$,
and the constraint vector $b \in \reals^m$.

Conic duality tells us that, if there exists any point in the interior
of $K$, \ie, $\intr K \ne \emptyset$, then this problem and the following
problem, called the \emph{dual problem}, have the same optimal value
\[
    \begin{aligned}
        & \text{maximize} && \bar f(A^Ty) + b^Ty\\
        & \text{subject to} && {-A^Ty} \in K^*,
    \end{aligned}
\]
with variable $y \in \reals^n$, where
\[
    \bar f(z) = \inf_x (f(x) - x^Tz),
\]
is sometimes known as the \emph{concave conjugate}. As we only use this fact
once in the main text, we do not derive it in detail, but
see~\cite[\S5.3.6]{bertsekas2009convex} for reference.

\section{Curve}\label{app:curve}
In this section, we derive the canonical trading function for a two-asset Curve
pool. Recall that the trading set for this market is given by~\cite{angeris2020improved}
\[
    S = \left\{R ~\middle|~ R_1 + R_2 - \frac{\alpha}{R_1 R_2} \ge k\right\}.
\]
From~\eqref{eq:canonical-reach}, we can write the trading function as
\[
    \phi(R) = \sup\left\{\lambda > 0 ~\middle|~ \frac{R_1 + R_2}{\lambda} - \frac{\alpha\lambda^2}{R_1 R_2} \ge k\right\}.
\]
Rewriting, we have that
\[
    \phi(R) = \sup\left\{\lambda > 0 ~\middle|~ 
    -\alpha\lambda^3 - k R_1 R_2 \lambda + R_1R_2(R_1 + R_2) \ge 0\right\}.
\]
The solution is given by the largest positive root of the cubic polynomial in
$\lambda$:
\[
    \lambda^\star = 
    \frac{\sqrt[3]{c_1(R) + \sqrt{c_2(R)}}}{3 \sqrt[3]{2} \alpha}
    - 
    \frac{\sqrt[3]{2} k R_1 R_2}{\sqrt[3]{c_1(R) + \sqrt{c_2(R)}}}
\]
where $c_1(R) = 27 \alpha^2 R_1^2 R_2 + 27 \alpha^2 R_1 R_2^2$ 
and $c_2(R) = 108 \alpha^3 k^3 R_1^3 R_2^3 + c_1^2$.
Plugging this back in, we have the canonical trading function
\[
    \phi(R) = \frac{R_1 + R_2}{k\lambda^\star} - \frac{\alpha(\lambda^\star)^2}{kR_1 R_2},
\]
which can (painfully) be verified to be homogeneous.

\section{Proof of concavity of Uniswap v3}\label{app:concavity}
The main difficulty in showing that~\eqref{eq:uni-hom} is concave is the square
root term
\[
    \sqrt{(\beta R_1 + \alpha R_2)^2 + 4(k-\alpha\beta)R_1R_2}.
\]
Its concavity follows from the fact that the set
\[
    Q = \{(x, y, t) \in \reals_+^3 \mid \|(\sqrt{\eta}(x-y), t)\|_2 \le \sqrt{1+\eta}(x+y) \}
\]
is convex when $\eta \ge 0$, where $\|\cdot\|_2$ denotes the Euclidean norm.
(To see this, note that norms are convex and affine functions are convex. Sets
of the form $\{z \mid f(z) \le 0\}$ are convex when $f$ is convex, and affine
precomposition preserves convexity.) Expanding the inequality gives the
following equivalent characterization of the set:
\[
    Q = \{(x, y, t) \in \reals_+^3 \mid t \le \sqrt{(x+y)^2 + 4\eta xy}\},
\]
which means that the function
\[
    \sqrt{(x+y)^2 + 4\eta xy} = \sup\{t \ge 0 \mid (x, y, t) \in Q\},
\]
is concave in $(x, y)$. Finally, setting $\eta = (k-\alpha\beta)/\alpha\beta$,
$x = \beta R_1$ and $y = \alpha R_2$ shows that the function
\[
    \sqrt{(\beta R_1 + \alpha R_2)^2 + 4(k - \alpha\beta)R_1R_2},
\]
is concave in $R_1$ and $R_2$.

\end{document}